\documentclass[a4paper,11pt]{amsart}

\usepackage[utf8]{inputenc}
\DeclareUnicodeCharacter{00A0}{ }
\usepackage{amssymb}
\usepackage{MnSymbol}
\usepackage{enumitem}

\usepackage{ae}
\usepackage[utf8]{inputenc}
\numberwithin{figure}{section}

\usepackage{mathrsfs,a4wide}
\usepackage{esint}
\usepackage{color}

\usepackage{comment}

\usepackage[colorlinks=true, linkcolor=red, urlcolor=red,citecolor=red]{hyperref}

\usepackage{color}

\usepackage{import}

\usepackage[leqno]{amsmath}

\numberwithin{figure}{section}

\newtheorem{theorem}{Theorem}[section]

\newtheorem{lemma}[theorem]{Lemma}
\newtheorem{proposition}[theorem]{Proposition}

\theoremstyle{definition}

\newtheorem{remark}[theorem]{Remark}
\numberwithin{equation}{section}

\addtocontents{toc}{\setcounter{tocdepth}{1}}

\newcommand{\D}{\mathrm{D}}
\newcommand{\R}{\mathbb{R}}
\newcommand{\N}{\mathbb{N}}

\newcommand{\Ha}{\mathcal{H}}

\newcommand{\beq}{\begin{equation}}
\newcommand{\eeq}{\end{equation}}
\newcommand{\dist}{{\rm dist}}
\newcommand{\eps}{\varepsilon}
 
\newcommand{\la}{\langle}
\newcommand{\ra}{\rangle}
\newcommand{\dia}{\mathrm{diam}}

\newcommand{\Div}{\operatorname{div}}
\newcommand{\pa}{\partial}

\newcommand{\spt}{\mathrm{spt}}
\newcommand{\medint}{-\kern -,375cm\int}
\newcommand{\medintinrigo}{-\kern -,315cm\int}

\renewcommand{\d}{\mathrm{d}}

\begin{document}

\title[Quantitative Alexandrov theorem]{Quantitative Alexandrov theorem and asymptotic behavior of the volume preserving mean curvature flow}

\author{Vesa Julin}

\author{Joonas Niinikoski}


\begin{abstract} 
We prove a new quantitative version of the Alexandrov theorem which states that if the mean curvature of a regular set in $\R^{n+1}$ is close to a constant in $L^{n}$-sense, then the set is close to a union of disjoint balls with respect to the  Hausdorff distance.   This result is more general than the previous quantifications of the Alexandrov theorem  and   using it we are able to  show that in $\R^2$ and   $\R^3$ a weak solution of the  volume preserving mean curvature flow starting from a set of finite perimeter  asymptotically convergences to a disjoint union of equisize balls, up to  possible translations.    Here by  weak solution we mean a flat flow, obtained via the minimizing movements scheme.  
\end{abstract}

\maketitle

\tableofcontents


\section{Introduction}

The main purpose of this article is to study the asymptotic behavior of  the weak solution of  the volume preserving mean curvature flow starting from  a set of finite perimeter. In the classical setting we are given a smooth set $E_0 \subset \R^{n+1}$ and we let  it evolve into a smooth family of sets  $(E_t)_t$ according to the law, where the normal velocity $V_t$ is proportional to the mean curvature of $E_t$ as
\beq
\label{the flow}
V_t = -(H_{E_t} - \bar H_{E_t}) \qquad \text{on }\, \pa E_t, 
\eeq
where $\bar H_{E_t} = \fint_{\pa E_t} H_{E_t} \, d \Ha^n$. Mean curvature type of equations are important in geometry, where one usually  studies the geometric properties   of $\pa E_t$ which are inherited from $\pa E_0$. The equation \eqref{the flow}  can also be seen as a volume preserving   gradient flow of the surface area. These equations arise naturally in physical models involving surface tension (see \cite{TCH}).

The main issue with \eqref{the flow} is that it may develop singularities in finite time even in the plane \cite{May, MS}. In order to pass over the singular time one may try to do a surgery procedure and restart the flow after a singular time as in \cite{HS} or to define a weak solution of \eqref{the flow}, which is what we will consider here. For the mean curvature flow one may define a weak solution by using the varifold setting by Brakke   \cite{Bra}, the level set  solution developed independently by  Chen-Giga-Goto \cite{CGG} and Evans-Spruck \cite{ES}, or by using the minimizing movements scheme  developed independently by Almgren-Taylor-Wang \cite{ATW} and  Luckhaus-St\"urzenhecker \cite{LS}. Since we want the solution of  \eqref{the flow} to be a family of sets and since \eqref{the flow} does not satisfy the comparison principle, the natural  choice is   to define a weak solution via the minimizing movements scheme as in  \cite{ATW, LS}. This solution is usually called a flat flow and it is well-defined due to \cite{MSS}, but not in general unique.

The advantage of the flat flow is that it is defined for all times for any bounded initial set with  finite perimeter and we may thus study its asymptotic behavior. Heuristically, one may guess that  the flat flow converges to a critical point of the static problem, which are classified in the recent work by Delgadino-Maggi \cite{DM} as disjoint union of balls, possibly tangent to each other. The asymptotic convergence  of \eqref{the flow} has been proved for initial sets with certain geometric properties such as convexity \cite{Hu}, nearly spherical  \cite{EsS} or sets which are near a stable critical set in the flat torus in low dimensions \cite{Joonas}.  We note that in these cases the flow does not develop singularities and is thus classically well-defined for all times. The result in  \cite{KK} shows that the convergence holds also for star-shaped sets, up to possible translations.    For the mean curvature flow with  forcing  the asymptotic behavior has been studied for the level set solution in \cite{GMT, GTZ}	
and for the flat flow in the plane in \cite{FJM}. The result closest to ours is  the recent work by Morini-Ponsiglione-Spadaro  \cite{MoPoSpa}, where the authors prove that the discrete-in-time  approximation of the flat flow  of \eqref{the flow} converges exponentially fast to disjoint union balls.  Here we are able to pass the time discretization to zero and characterize the limit sets for the flat flow of \eqref{the flow} in  $\R^2$ and $\R^3$. The precise definition of the flat flow is given in Section 4.  
\begin{theorem}
\label{thm2}
Assume $E_0 \subset \R^{n+1}$, with $n\leq 2$,  is a bounded set of finite perimeter which is either essentially open or essentially closed and let $(E_t)_{t \geq 0}$ be a flat flow of \eqref{the flow} starting from $E_0$.  There is $N \in \N$ 
such that the following holds: for every $\eps >0$ there is $T_\eps>0$ such that  for every $t \geq T_\eps$ there are points $x_1, \dots, x_N$, which may depend on time,  with  $|x_i -x_j| \geq  2r$  for $i \neq j$ and 
$r=N^{-\frac{1}{n+1}}$ such that   for $F_t = \bigcup_{i=1}^N B_r(x_i)$ it  holds 
\[
\sup_{x \in \pa E_t} d_{\pa F_t}(x) \leq \eps .
\]
\end{theorem}
Here $d_{\pa F}$ denotes the distance function.  To the best of our knowledge this is the first result on the characterization of  the asymptotic limit  of \eqref{the flow} in $\R^3$. The above result holds for any limit of the approximative flat flow and we  do not need the additional assumption on the convergence of the perimeters as in \cite{LS, MSS}.  We note that the assumption on $E_0$ being either  essentially open or closed is only needed to ensure that the flow is continuous up to time zero. It plays no role in the asymptotic analysis. 

Concerning the limiting configurations,   Theorem \ref{thm2} is sharp since  the flow \eqref{the flow} may converge to tangent balls as  it is shown in \cite{FJM}. On the other hand, we believe that one can rule out the possible translations and the flow  actually convergences to a disjoint union of balls. The higher dimensional    case and the possible speed of convergence are also open problems.  

\subsection*{Quantitative Alexandrov theorem}

The proof of Theorem \ref{thm2} is based on the  dissipation inequality  proven in \cite{MSS} and  stated  in  Proposition \ref{spadaro0}. This implies that there is a sequence of times $t_j \to \infty$ such that  the mean curvatures of the evolving sets $E_{t_j}$ are asymptotically close to a constant with respect to the $L^2$-norm. Therefore, we need a quantified version of the Alexandrov theorem which  enables us to conclude that the sets $E_{t_j}$ are close to a disjoint union  of  balls.

There is  a lot of recent  research on generalization of the Alexandrov theorem  \cite{CM, DM, DMMN,  RKS,  KM}. We refer  the survey paper \cite{Ci} for the state-of-the-art. Unfortunately, none of the available results is applicable to our problem, and we are also not able to use the characterization of the critical sets by Delgadino-Maggi \cite[Corollary 2]{DM} to identify  the limit set. Indeed, even if we know that the sets $E_{t_j}$ converge to a set of finite perimeter and their mean curvatures converge to a constant, it is not clear why the limit set is a set of finite perimeter with weak  mean curvature as this class of sets is not in general closed. Our main result of the paper is the following quantification of the Alexandrov theorem, which is the main technical tool in the proof of  Theorem \ref{thm2}.     
\begin{theorem}
\label{thm1}
Let $E \subset \R^{n+1}$ be a $C^2$-regular set such that $P(E) \leq C_0$ and $|E| \geq 1/C_0$. There are positive constants $q=q(n) \in (0,1] $, $C = C(C_0,n)$ and $\delta = \delta(C_0,n)$
such that if $\| H_E - \lambda \|_{L^{n}(\pa E)} \leq \delta$ for some $\lambda \in \R$, then 
$1/C \leq \lambda \leq C$ and there are points  $x_1, \dots, x_N$  with $|x_i - x_j| \geq 2R$, where $R=n/\lambda$, such that for $F = \bigcup_{i=1}^N B_{R}(x_i)$ it holds
\[
\sup_{x \in \pa E} d_{\pa F}(x) \leq C  \| H_E - \lambda \|_{L^{n}(\pa E)}^q.
\]
 Moreover,
\[
\Big{|} P(E) -  N (n+1) \omega_{n+1} R^{n}  \Big{|} \leq C \| H_E - \lambda \|_{L^{n}(\pa E)}^q.
\]
\end{theorem}
The main advantage  of Theorem \ref{thm1}  with respect to the previous results in the literature is that we do not assume any geometric restriction on $E$  such as mean convexity. Moreover, we
assume the mean curvature to be close to a constant only in the $L^n$-sense, which is exactly what we need for the asymptotic analysis in Theorem \ref{thm2}. This makes the proof challenging
as we, e.g.,  cannot use the estimates from the Allard regularity theory \cite{All}.

Theorem \ref{thm1} is sharp in the sense that $ \| H_E - \lambda \|_{L^{n}(\pa E)}$ cannot be replaced by a weaker $L^p$-norm.  This can be easily seen by considering a  set which is  a union of the unit ball and a  ball of small radius $\eps$ far away. On the other hand, the dissipation inequality in  Proposition \ref{spadaro0} controls only the $L^2$-norm of the mean curvature, which is the reason why we cannot prove Theorem \ref{thm2} in higher dimensions. The proof of Theorem \ref{thm1} is done in a constructive way and we obtain an explicit bound on the  exponent $q= (n+2)^{-3}$. It would be interesting to obtain the sharp one as it might be crucial in order to obtain the possible exponential convergence of \eqref{the flow}  as    in \cite{MoPoSpa}.  In the two-dimensional case the optimal power $q=1$ is proven in \cite{FJM}.

\subsection*{Outline of the proof of Theorem \ref{thm1}}

Since the proof of Theorem \ref{thm1} is rather long, we give its  outline here. The argument is based on the proof of the Heinze-Karcher inequality by Montiel-Ros \cite{MR}, which is originally  an alternative proof for \cite{Ros}.  We revisit the argument by Montiel-Ros and deduce in Proposition \ref{monros}  that for $E$ and $R$ as in Theorem \ref{thm1} and for  $0 < r <R$  the volume of the  set $E_r = \{ x \in E : \text{dist}(x, \pa E) >r \}$ satisfies the estimate  
\[
\Big{|}  |E_r|  - \frac{|E|}{R^{n+1}} (R-r )^{n+1}  \Big{|} \leq  C \| H_E - \lambda\|_{L^n(\pa E)}. 
\]
We use this in Step 1 of the proof of Theorem \ref{thm1} to deduce that for $r$ close to $R$ the set $E_r$ is a union of finite number of components, or clusters, with positive distance to each other. 

We note that the above inequality is not enough to conclude the proof as, for example,    the  cube $Q = (-1,1)^{n+1}$ satisfies $|Q_r| = (1-r)^{n+1}|Q|$. Therefore, we need further information from the Montiel-Ros argument and we prove in Proposition \ref{monros} that  the Minkowski sum  $E_r + B_\rho = \{x \in \R^{n+1} : \text{dist}(x,E_r) < \rho \}$,  with $0<\rho<r<R$, satisfies 
\[
\Big{|}  |E_r + B_\rho|  - \frac{|E|}{R^{n+1}} (R-(r-\rho) )^{n+1}  \Big{|} \leq   \frac{C}{(R-r)^{n+1}}  \| H_E - \lambda\|_{L^n(\pa E)}.
\]
This enables us to prove that the components of $E_r + B_\rho \subset E$, with properly chosen  $\rho$ and $r$, are almost spherical. In particular, if $E$ satisfies the above estimate with $C=0$, then it is a disjoint union of balls. This, together with the density estimate   from  \cite{Top}, concludes the proof.


\section{Notation and preliminary results}

In this section we briefly introduce  our notation and recall some  results  from differential geometry. Given a set $E \subset \R^{n+1}$   the distance function $d_E :  \R^{n+1}  \to [0,\infty)$ is defined, as usual, as
\[
d_E(x) := \inf_{y \in E} |x - y|
\] 
and we denote the signed distance function by $\bar{d}_E : \R^{n+1} \to \R$, 
\[
\bar{d}_E(x) := \begin{cases} - d_{\pa E}(x) , \,\,&\text{for }\, x \in E \\
d_{ \pa E}(x) , \,\, &\text{for }\, x \in \R^{n+1} \setminus E.  \end{cases} 
\]
Then clearly it holds $d_{\pa E} = |\bar d_E|$. We denote the ball with radius $r$ centered at $x$ by $B_r(x)$ and by $B_r$ if it is centered at the origin. Given a set $E \subset \R^{n+1}$ we denote  its $\rho$-enlargement by the Minkowski sum  
\[
E + B_\rho = \{ x + y \in \R^{n+1} : x \in E , \, \,  y \in B_\rho \} = \{x \in \R^{n+1} : d_E(x) < \rho \}.
\] 

For a measurable set $E \subset \R^{n+1}$ the shorthand notation $|E|$ denotes its Lebesgue measure and we  denote  the $k$-dimensional measure of the unit ball in $\R^k$  by $\omega_k$. In some cases, we may use the shorthand notation $|E|$ more generally for a measurable set $E \subset \R^k$ to denote its $k$-dimensional  Lebesgue measure but this shall be clear from context.

For a set of finite perimeter  $E \subset \R^{n+1}$  we denote its reduced boundary by $\pa^* E$ and the perimeter by $P(E)$. Recall that $P(E) = \Ha^{n}(\pa^* E)$ and  for regular enough set it holds $\pa^* E = \pa E$. The relative isoperimetric inequality states that for every set of finite perimeter $E$ and for every ball $B_r(x)$ it holds
\[
\Ha^n(\pa^* E \cap B_r(x))^{\frac{n+1}{n}} \geq c_n \, \min \big{\{} |E \cap B_r(x)|,  |B_r(x) \setminus E| \big{\}},
\]
for a dimensional constant.  We refer to  \cite{Ma} for an introduction to the topic.   

 We  define  the \emph{tangential differential}  of $F \in C^1(\R^{n+1}; \R^m)$  on $\pa E$ by 
\[
D_\tau F(x) = D F(x) (I - \nu_E(x)\otimes \nu_E(x)) ,
\]
where $\nu_E$  denotes the unit outer normal of $E$. For a function $f \in C^1(\R^{n+1}; \R)$ we denote by $\nabla_\tau f$ its tangential gradient which is a vector in $\R^{n+1}$.   We define the tangential divergence of $F \in C^1(\R^{n+1}; \R^{n+1})$ by $\text{div}_\tau F = \text{Tr} (D_\tau F)$ . Then the divergence theorem on manifolds  generalizes to  
\[
\int_{\pa^* E} \text{div}_\tau F  \, d \Ha^n = \int_{\pa^* E} H_E \,  \la F, \nu_E \ra  \, d \Ha^n,
\]
where $H_E \in L^1(\pa^* E)$ is the  distributional mean curvature. When $\pa E$ is smooth $H_E$ agrees with the classical definition of the mean curvature,  which for us   is the sum of the principal  curvatures. 

We begin by recalling the well-known inequality proven first by  Simon \cite{Sim} in $\R^3$ and then  by  Topping \cite{Top} in the general case. 

\begin{theorem}
\label{Topping}
Let \(\Sigma \subset \R^{n+1}\) be a compact and connected \(C^2\)-hypersurface. Then
\beq
\label{topping}
\dia (\Sigma) \leq C_n \int_\Sigma|H_\Sigma|^{n-1} \ \d \Ha^n,
\eeq
where \(C_n\) depends only on the dimension.
\end{theorem}

We need also the Michael-Simon inequality \cite{MiSi}.

\begin{theorem}
\label{M-S} Let \(\Sigma \subset \R^{n+1}\), \(n\geq 2\), be a compact  \(C^2\)-hypersurface. Then 
for every non-negative 
\(\varphi \in C^1(\mathbb R^{n+1})\) 
\beq
\label{ms}
\|\varphi\|_{L^\frac{n}{n-1}(\Sigma)} \leq C_n \int_\Sigma |\nabla_\tau \varphi| + \varphi |H_\Sigma| \ \d \Ha^n,
\eeq
where \(C_n\) depends only on the dimension.
\end{theorem}

The following  density-type estimate is essentially proven in \cite[Lemma 2.1]{MoPoSpa}.

\begin{proposition}
\label{scalecontrol}
Let \(E \subset \mathbb R^{n+1}\) be a set of finite perimeter with \(P(E)>0\) and \(0<\beta<1\). There is a positive constant \(c=c(n,\beta)\) such that
\[
r_{E,\beta}:= \sup\left\{r \in \mathbb R_+ : \text{there is \(x \in \mathbb R^{n+1}\) with} \ |B_r(x) \cap E| \geq \beta |B_r(x)|\right\} \geq c \frac{|E|}{P(E)}.
\]
\end{proposition}

We use the previous results to prove the following lemma, which is useful when we bound the Lagrange multipliers  and the number of the components of the flat flow of \eqref{the flow}. 

\begin{lemma}
\label{lambda-diam control}
Let $E \subset \R^{n+1}$ be a bounded set of finite perimeter with a distributional mean curvature $H_E \in L^1(\pa^* E)$,
$\lambda \in \R$ and $1 \leq C_0 < \infty$. There is a positive constant $C=C(C_0,n)$ 
such that the following hold.
\begin{itemize}
\item[(i)]
 If $P(E) \leq C_0$ and $|E| \geq 1/C_0$, then
\[
1/C-C \|H_E-\lambda\|_{L^1(\pa^* E)} \leq \lambda \leq C+ C\|H_E-\lambda\|_{L^1(\pa^* E)}.
\]
\item[(ii)]
If $P(E) \leq C_0$, $|E| \geq 1/C_0$ and $E$ is $C^2$-regular, then the 
number of the components of $E$  is bounded 
by $C(1+\|H_E-\lambda\|_{L^n(\pa E)}^n)$ and their diameters are bounded by
 $C(1+\|H_E-\lambda\|_{L^{n-1}(\pa E)}^{n-1})$.
\end{itemize}
\end{lemma}

\begin{proof}
Our standing assumptions throughout the proof are  $P(E) \leq C_0$ and  $|E| \geq 1/C_0$. The perimeter bound and the global isoperimetric inequality yield
\[
|E| \leq c_n P(E)^{\frac{n+1}{n}} \leq c_n C_0^{\frac{n+1}{n}}.
\]

By the assumptions on $E$ and by the divergence theorems we compute for any vector field ${F \in C^1(\R^{n+1};\R^{n+1})}$
\beq 
\label{aux1}
\begin{split}
\lambda \int_E \Div  F \ \d x  &=  \int_{\pa^* E} \lambda \la F, \nu_E \ra \, \d\Ha^n \\
&=  \int_{\pa^* E} H_E \la F, \nu_E \ra \, \d\Ha^n +  \int_{\pa^* E} (\lambda - H_E) \la F, \nu_E \ra \, \d\Ha^n \\
&=  \int_{\pa^* E} \Div_\tau F \, \d\Ha^n +  \int_{\pa^* E} (\lambda - H_E) \la F, \nu_E \ra \, \d\Ha^n.
\end{split}
\eeq
Our goal is to construct a suitable vector field $F$ to obtain (i) from \eqref{aux1}. To this aim, 
we use first the isoperimetric inquality, Proposition \ref{scalecontrol} and a suitable continuity argument to find
positive $r_0=r_0(C_0,n)$, $R_0=R_0(C_0,n)$ and $r$ such that $r_0 \leq r \leq R_0$ and, by possibly  translating the coordinates, 
$|B_r \cap E| = |B_r|/2$. Again, it follows from the relative isoperimetric inequality that  $\Ha^n(\pa^* E \cap B_r) \geq c$ with some positive $c=c(C_0,n)$. 
Choose a decreasing $C^1$-function 
$f:\mathbb R \rightarrow \mathbb R$ for which
\[
f(t)=
\begin{cases}
(2r)^{-1}, & \ \text{for} \  t \leq \frac{3}{2} r \\
t^{-1}, & \ \text{for} \  t \geq  \frac{5}{2} r
\end{cases}
\]
and the conditions $f(t) \leq \min\{(2r)^{-1},t^{-1}\}$ , $|f'(t)| \leq (2r)^{-2}$ hold on $[\frac{3}{2}r,\frac{5}{2}r]$.
We define $F:\R^{n+1} \rightarrow \R^{n+1}$ by setting $F(x)=f(|x|)x$. Then  $F$ is a $C^1$-vector field with 
\begin{align*}
\D F (x) &= f(|x|) I + \frac{f'(|x|)}{|x|} \ x \otimes x, \qquad \text{for every} \  \ x \in \R^{n+1}, \\
\Div F (x) &= (n+1)f(|x|)  + f'(|x|)|x|, \qquad \text{for every} \  \ x \in \R^{n+1} \ \ \text{and}\\
\Div_\tau F(x) &= nf(|x|) + f'(|x|)\left(|x| - \frac{\langle x, \nu_{E}\rangle^2}{|x|} \right), \qquad  \text{for every} \  \ x \in \pa^* E.
\end{align*}
Then  $0<\Div  F \leq (n+1) (2r)^{-1}$ everywhere and  $\Div F =(n+1) (2r)^{-1}$ in $B_r$
so by using these and the earlier observations  we obtain
\beq
\label{aux2}
 \frac{n+1}{4R_0} |B_{r_0}| \leq \frac{n+1}{4r} |B_r| \leq \frac{n+1}{2r} |B_r \cap E| \leq \int_E \Div  F \ \d x \leq \frac{n+1}{2r}|E|
\leq \frac{c_n (n+1)}{2r_0}C_0^\frac{n+1}{n}.
\eeq
 Again, $ 0 \leq \Div_\tau F \leq  n(2r)^{-1}$ on $\pa^* E$ and  $\Div_\tau F =  n(2r)^{-1}$ on $\pa^* E \cap B_r$ and thus
\beq 
\label{aux3}
\frac{nc}{2R_0} \leq \frac{n}{2r}\Ha^n(\pa^* E \cap B_r) \leq \int_{\pa^* E} \Div_\tau F \, \d\Ha^n \leq \frac{nP(E)}{2r} \leq \frac{nC_0}{2r_0}.
\eeq
We use  \eqref{aux1}, \eqref{aux2}, \eqref{aux3}  and   $|F| \leq 1$ to  obtain (i).

The claim (ii) is easy to prove in the planar case and therefore we assume that $n\geq 2$. Let  $E_1, E_2, \dots,E_N$ denote the connected components of $E$. 
 We apply Theorem \ref{M-S} on \(\pa E_i\) 
with \(\varphi=1\) and use Hölder's inequality to obtain 
\[
C_n^{-1} \leq \| H_{E_i} \|_{L^n(\pa E_i)} \leq \| H_{E_i} - \lambda \|_{L^n(\pa E_i)} +|\lambda| P(E_i)^\frac{1}{n},
\]
from which we conclude using (i) and Hölder's inequality 
\beq
\begin{split}
\label{aux4}
N C_n^{-n}  
&\leq 2^n\| H_E - \lambda \|_{L^n(\pa E)}^n +2^n|\lambda|^n P(E) \\
&\leq 2^n \| H_E - \lambda \|_{L^n(\pa E)}^n +2^{2n}C_0C^n\left(1+\| H_E - \lambda \|_{L^1(\pa E)}^n\right)  \\
&\leq 2^n \| H_E - \lambda \|_{L^n(\pa E)}^n +2^{2n}C_0C^n\left(1+C_0^{n-1}\| H_E - \lambda \|_{L^n(\pa E)}^n\right) .
\end{split}
\eeq
On the other hand, Theorem \ref{Topping} together with (i) and Hölder's inequality implies 
\beq
\begin{split}
\label{aux5}
\sum_{i} \text{diam}(E_i) &\leq  \sum_{i} C_n \int_{\pa E_i}|H_{E_i}|^{n-1} \, \d \Ha^n \\
&\leq  \sum_{i}  2^{n-1}C_n\left(\int_{\pa E_i}|H_{E_i}-\lambda|^{n-1} \, \d\Ha^n + |\lambda|^{n-1} P(E_i)\right) \\
&\leq 2^{n-1}C_n\left(\int_{\pa E}|H_E-\lambda|^{n-1} \, \d\Ha^n + P(E) |\lambda|^{n-1}\right) \\
&\leq 2^{n-1}C_n\left(\| H_E - \lambda \|_{L^{n-1}(\pa E)}^{n-1}+ 2^{n-1}C_0C^n(1+\| H_E - \lambda \|_{L^1(\pa E)}^{n-1})\right) \\
&\leq 2^{n-1}C_n\left(\| H_E - \lambda \|_{L^{n-1}(\pa E)}^{n-1}+ 2^{n-1}C_0C^n(1+C_0^{n-2}\| H_E - \lambda \|_{L^{n-1}(\pa E)}^{n-1})\right).
\end{split}
\eeq
Thus, by possibly increasing $C$, the latter part of the claim follows from \eqref{aux4} and  \eqref{aux5}.
\end{proof}


\section{Quantitative Alexandrov theorem }

We split the proof of  Theorem \ref{thm1} into two parts. We first revisit the Montiel-Ros argument in Proposition \ref{monros} where all the technical heavy lifting is done. The idea of Proposition \ref{monros}  is to transform the (local) information of  the mean curvature of  $E$ being close to a constant, into information on the $\rho$-enlargement of the level sets  of the distance function of $\pa E$. We note that the statement of  Proposition \ref{monros} is given by the sharp exponent. The proof of Theorem \ref{thm1} is then  based on purely   geometric arguments. 

We first state the following equivalent formulation of the theorem. 

\begin{remark}
\label{rem thm1}
Once we prove that in Theorem \ref{thm1} the number of component of $E$ is bounded, the statement on the $L^\infty$-distance is equivalent to the fact that,  
under the assumption $\| H_E - \lambda \|_{L^{n}(\pa E)} \leq \delta$, there are points  $x_1, \dots, x_N$ such that
\[
\bigcup_{i=1}^N B_{\rho_-}(x_i) \subset E \subset \bigcup_{i=1}^N B_{\rho_+}(x_i),
\]
where $\rho_-=R-C\| H_E - \lambda \|_{L^{n}(\pa E)}^q$, $\rho_+=R+C\| H_E - \lambda \|_{L^{n}(\pa E)}^q$, $R=n/\lambda$ and the balls 
$ {B_{\rho_-}(x_1),\ldots, B_{\rho_-}(x_N)}$ are disjoint to each other.    We leave the details  to the reader.
\end{remark}

In Theorem \ref{thm1} we  assume that the mean curvature is bounded only in the $L^n$-sense and thus the estimates from the Allard's regularity theory \cite{All} are not available for us. Indeed,  the $L^n$-boundedness of the mean curvature  is not strong enough to give  proper density estimates. Moreover, even in the three dimensional case  $\R^3$ we cannot use the results from \cite{Sim}, because we do not have a uniform bound on the Euler characteristic of the set $E$.    However, if we know that the mean curvature is close to a constant with respect to the $L^n$-norm, then  the following density estimate holds.  The proof is based on \cite[Lemma 1.2]{Top}. 
\begin{lemma}
\label{density}
Let \(\Sigma \subset \R^{n+1}\) be a compact  \(C^2\)-hypersurface and \(\lambda \in \mathbb R_+\). 
There is a positive dimensional constant $\delta_n$  such that if \(\|H_\Sigma-\lambda\|_{L^n(\Sigma)} \leq \delta_n\), then 
\[
\delta_n \leq \frac{\Ha^n(B(x,r) \cap \Sigma)}{r^n}
\]
for every \(x \in \Sigma\) and \(0<r \leq \frac{\delta_n}{\lambda}\).
\end{lemma}

\begin{proof}
The planar case $n=1$ is rather obvious  and we leave it to the reader.
Let us assume $n\geq2$.  Fix \(x \in \Sigma\) and define \(V:[0,\infty) \rightarrow [0,\infty)\) as \(V(r)=\Ha^n(B_r(x) \cap \Sigma)\).
Since \(V\) is increasing, the derivative \(V'(r)\) is defined for almost every \(r \in [0,\infty)\) and
\[
\int_{r_1}^{r_2} V'(\rho) \ \d \rho \leq V(r_2) - V_2(r_1) \ \ \text{whenever} \ \ 0 \leq r_1 < r_2.
\]
By a standard foliation argument we have that \(\Ha^n(\partial B_r(x) \cap \Sigma)>0\) at most countably many 
\(r \in \mathbb R_+\). Thus \(V'(r)\) is defined and \(\Ha^n(\partial B_r(x) \cap \Sigma)=0\) for almost every \(r \in [0,\infty)\).
Fix such \(r\) and choose \(h \in \mathbb R_+\) for which \(\Ha^n(\partial B_{r+h}(x) \cap \Sigma)=0\). Define a cut-off function
\(f_h: \mathbb R^{n+1} \rightarrow \mathbb R\) by setting
\[
f_h(y)=
\begin{cases}
1, &y \in B_r(x) \\
1-\frac{|y-x|}{h}, &y \in B_{r+h}(x)\setminus B_r(x) \\
0, &y \in \mathbb R^{n+1} \setminus B_{r+h}(x).
\end{cases}
\]
By using a suitable approximation argument combined with Theorem \ref{M-S} we obtain
\[
V(r)^\frac{n-1}{n} \leq  C_n \left( \frac{V(r+h)-V(r)}{h} + \|f_h H_\Sigma \|_{L^1(\Sigma)}\right).
\]

In turn, we may choose a sequence 
\((h_k)_{k} \) such that \(h_k \rightarrow 0\) and \(\Ha^n (\partial B_{r+h_k}(x) \cap \Sigma)=0\). Then by 
letting \(k \rightarrow \infty\) the previous estimate yields
\begin{align*}
V(r)^\frac{n-1}{n} 
&\leq C_n \left( V'(r) + \int_{\overline B_r(x) \cap \Sigma}|H_\Sigma| \ \d \Ha^n \right) \\
& \leq C_n \left(V'(r) + \int_{B_r(x) \cap \Sigma}|H_\Sigma| \ \d \Ha^n\right) \\
& \leq C_n \left(V'(r) + \int_{B_r(x) \cap \Sigma}|H_\Sigma - \lambda| \ \d \Ha^n + \lambda V(r)\right) \\
& \leq C_n \left(V'(r) +  \|H_\Sigma - \lambda\|_{L^n(\Sigma)} V(r)^\frac{n-1}{n} + \lambda V(r)\right).
\end{align*}
 Thus for almost every \(r \in   (0,\infty)\) it holds
\[
\left(\frac{C_n^{-1}- \|H_\Sigma - \lambda\|_{L^n(\Sigma)}}{V(r)^\frac{1}{n}} - \lambda \right) V(r) \leq V'(r).
\]
If $\|H_\Sigma - \lambda\|_{L^n(\Sigma)} \leq \delta_n$  for  small $\delta_n$ then the above inequality implies
\[
 \frac{1}{2C_n} V(r)^{1-\frac{1}{n}}- \lambda  V(r) \leq V'(r).
\]

Fix $r <\delta_n/\lambda$. We  assume that  $V(r) \leq \delta_n r^n$, since otherwise the claim is trivially true. By the monotonicity we have $V(\rho)^\frac1n  \leq V(r)^\frac1n  \leq \delta_n/\lambda $ for all $0 < \rho < r$.  For $\delta_n$  small enough the above   inequality then yields
\[
 \frac{1}{4 C_n}V(\rho)^{1-\frac{1}{n}} \leq V'(\rho)
\]
for almost every $0 < \rho < r$. The claim follows by integrating this over $(0,r)$. 

\end{proof}

\subsection{Montiel-Ros argument}

We  recall that for $E \subset \R^{n+1}$ we denote
\beq \label{E_r}
E_r := \{ x \in E : \text{dist}(x, \pa E) > r \}.
\eeq
We use the fact that  $E$ is $C^2$-regular and say that $x \in \pa E$ satisfies interior ball condition with radius $r$, if for  $y = x - r\nu_E(x)$ it holds $B_r(y) \subset E$. For $r>0 $ we  define
\beq \label{gamma_r}
\Gamma_r := \{ x \in \pa E : x \,\, \text{satisfies interior ball condition with radius } \,\, r \}.
\eeq

\begin{proposition}
\label{monros}
Let $\lambda \in \R$ and suppose that a  bounded and  \(C^2\)-regular set  $E \subset \R^{n+1}$  satisfies $P(E) \leq C_0$ and $|E| \geq 1/C_0$ with $C_0 \in \R_+$.
Then for $0< r < R$  with $R = n/\lambda$ it holds
\[
\Big{|}  |E_r|  - \frac{|E|}{R^{n+1}} (R-r )^{n+1}  \Big{|} \leq  C \| H_E - \lambda\|_{L^n(\pa E)}
\]
and 
\[
\Ha^{n}(\pa E \setminus \Gamma_{r}) \leq \frac{C}{(R-r)^{n+1}} \| H_E - \lambda\|_{L^n(\pa E)},
\]
provided that $\| H_E - \lambda\|_{L^n(\pa E)} \leq \delta$, where the constants $C$ and $\delta$ depend only on $C_0$ and on the dimension. Moreover, under the same assumptions,  for $0< \rho < r < R$   it holds
\[
\Big{|}  |E_r + B_\rho|  - \frac{|E|}{R^{n+1}} (R-(r-\rho) )^{n+1}  \Big{|} \leq   \frac{C}{(R-r)^{n+1}}  \| H_E - \lambda\|_{L^n(\pa E)}.
\]
\end{proposition}

\begin{proof} 
As we already mentioned the proof is based on   the Montiel-Ros argument for the Heinze-Karcher inequality, which we recall shortly. To that aim we define $\zeta : \pa E \times \R \to \R^{n+1}$ as 
\[
\zeta(x,t) = x - t \nu_E(x).
\]
 We denote the principle curvatures of $\pa E$ at $x$ by $k_1(x), \dots k_n(x)$ and assume that they are pointwise ordered as $k_i(x) \leq k_{i+1}(x)$.  If we consider $\pa E \times \R$ as a hypersurface embedded in $\R^{n+2}$ then its tangential Jacobian is 
\[
J_\tau \zeta(x,t) = \prod_{i=1}^n |1 - tk_i(x)| \qquad \text{on } \, \pa E \times \R.
 \] 
For every bounded Borel set $M \subset \pa E \times \R$ we have by the area formula
\[
\int_{\zeta(M)} \Ha^0(\zeta^{-1}(y)\cap M)\, \d y = \int_M J_\tau \zeta \, \d\Ha^n.
\]
In the proof, \(C\) denotes a positive constant which may change from line to line,  depending only on \(C_0\) and on the dimension.
\ \\
\ \\
\textbf{Step 1:} \quad 
In order to utilize Lemma \ref{lambda-diam control}, we choose $\delta=\delta(C_0,n)$ to be same as in
the lemma and assume $\| H_E - \lambda\|_{L^n(\pa E)} \leq \delta$.
Then  $E$ has $N$ many connected components with $N \leq C$. We may thus prove the claim componentwise and   assume that $E$ is connected.
We denote 
\[
\Sigma := \{ x \in \pa E : | H_E(x) - \lambda | < \lambda/2\}. 
\]
By Lemma \ref{lambda-diam control} it holds $\lambda \geq 1/C$ and thus  by H\"older's inequality  it  holds 
\beq 
\label{monros1}
\Ha^n(\pa E \setminus \Sigma) \leq \frac{2}{\lambda} \int_{\pa E}| H_E(x) - \lambda |  \, \d \Ha^n \leq C \| H_E - \lambda\|_{L^n(\pa E)}.
\eeq
Moreover, we have 
\[
\begin{split}
\frac{n}{n+1} \int_{\Sigma} \frac{1}{H_E}\, \d \Ha^n &= \frac{n}{n+1} \int_{\Sigma}\left( \frac{1}{\lambda}  + \left(\frac{1}{H_E} - \frac{1}{\lambda} \right)\right) \, \d \Ha^n\\
&\leq \frac{n P(E)}{(n+1) \lambda} + C \| H_E(x) - \lambda\|_{L^n(\pa E)}.
\end{split}
\]
Since \(E\) is connected, Lemma \ref{lambda-diam control} yields $\text{diam}(E) \leq \tilde R$ with $\tilde R= \tilde R(C_0,n) \geq R$. Choose $x_0 \in E$. Then using  \eqref{aux1}  with $F(x) = x - x_0$ we obtain
\[
nP(E) =  (n+1)\lambda \,  |E| + \int_{\pa E} (H_{E} - \lambda)  \la (x-x_0), \nu_E \ra \, \d\Ha^{n} ,
\]
which in turn implies
\beq
\label{peri-vol}
\big{|}nP(E) - (n+1)\lambda \,  |E|  \big{|} \leq C \| H_E - \lambda\|_{L^n(\pa E)}.
\eeq
Hence, we deduce
\beq 
\label{monros2}
\frac{n}{n+1} \int_{\Sigma} \frac{1}{H_E}\, \d \Ha^n \leq  |E|  + C \| H_E - \lambda\|_{L^n(\pa E)}.
\eeq

Next we define
\[
Z = \{ (x,t) \in \Sigma \times [0,\infty) : 0 \leq t \leq 1/k_n(x) \}.
\]
Note that this is well-defined, since $x \in \Sigma$ implies  $k_n(x) \geq \frac{H_E(x)}{n} \geq \frac{\lambda}{2n} >0$.  We also set 
\[
\Sigma'_1 = \{ x \in \pa E \setminus \Sigma : k_n(x) \leq 1/\tilde R\} \qquad \text{and} \qquad \Sigma'_2 = \{ x \in \pa E \setminus \Sigma : k_n(x) > 1/\tilde R\},
\]
\[
Z_1' = \Sigma'_1 \times [0,\tilde R] \qquad \text{and} \qquad Z_2' =   \{ (x,t) \in \Sigma_2' \times [0,\infty) : 0 \leq t \leq 1/k_n(x) \}
\]
and finally  
\[
Z' = Z_1' \cup Z_2' .
\]
Then $Z$ and $Z'$ are disjoint and bounded Borel sets and  it holds $E \subset \zeta(Z \cup Z')$. To see this fix $y \in E$ and  let $x \in \pa E$ be such that $r = d_{\pa E}(y) = |x-y|$. Then we may write 
$y = x - r \nu_E(x)$ and by the maximum principle  $k_n(x) \leq 1/r$. Since $\text{diam}(E) \leq \tilde R$, then  $r \leq \tilde R$ and we conclude that $(x,r) \in Z \cup Z'$ and $y=\zeta(x,r).$

We now recall the Montiel-Ros argument. We use  the fact that $E \subset \zeta(Z \cup Z')$, the area formula, the arithmetic geometric inequality and the fact that for $x \in \Sigma$ it holds  $1/k_n(x) \leq n/H_E(x)$ to obtain 
\[
\begin{split}
|E| &\leq |\zeta (Z)| + |\zeta (Z')|  \leq \int_{\zeta(Z)} \Ha^0(\zeta^{-1}(y) \cap Z) \, \d y + |\zeta (Z')| \\
&= \int_Z J_\tau \zeta \, \d \Ha^n + |\zeta (Z')|  \\
&=  \int_\Sigma \int_0^{1/k_n(x)} \prod_{i=1}^n (1 - tk_i(x))  \, \d t \d \Ha^n  + |\zeta (Z')|  \\
&\leq  \int_\Sigma \int_0^{1/k_n(x)} \left(1 - \frac{t}{n} H_E(x)\right)^n  \, \d t \d \Ha^n   + |\zeta (Z')| \\
&\leq  \int_\Sigma \int_0^{n/H_E(x)} \left(1 - \frac{t}{n} H_E(x)\right)^n  \, \d t \d \Ha^n   + |\zeta (Z')|= \frac{n}{n+1} \int_{\Sigma} \frac{1}{H_E} \, \d \Ha^n + |\zeta (Z')|  . 
 \end{split}
\]
Next we quantify the previous four inequalities. To that aim we define the non-negative numbers $R_1, R_2, R_3$ and $ R_4$ as
\begin{align}
R_1 &= |\zeta (Z) \setminus E| \label{R1} \\
R_2 &=\int_{\zeta(Z)} |\Ha^0(\zeta^{-1}(y) \cap Z) -1| \, \d y \label{R2} \\
R_3 &= \int_\Sigma \int_0^{1/k_n(x)} \Big{|} \left(1 - \frac{t}{n} H_E(x)\right)^n -  \prod_{i=1}^n (1 - tk_i(x)) \Big{|} \, \d t \d \Ha^n  \label{R3} \\
R_4 &= \int_\Sigma \int_{1/k_n(x)}^{n/H_E(x)} \Big{|}1 - \frac{t}{n} H_E(x)\Big{|}^n \, \d t \d \Ha^n. \label{R4}
\end{align}
Then by repeating the Montiel-Ros argument we deduce that
\[
|E| \leq  \frac{n}{n+1} \int_{\Sigma} \frac{1}{H_E} \, \d \Ha^n + |\zeta (Z')|  - R_1  - R_2 - R_3 - R_4. 
\]
Therefore, by \eqref{monros2} it holds 
\[
R_1  + R_2 + R_3 + R_4 \leq |\zeta (Z')| + C \| H_E - \lambda\|_{L^n(\pa E)},
\]
where $R_i$ are  defined in \eqref{R1}-\eqref{R4}. 


Let us next show that 
\beq
\label{Z'}
|\zeta (Z')| \leq C \| H_E(x) - \lambda\|_{L^n(\pa E)}.
\eeq
Indeed, by the area formula we have 
\beq \label{monros3}
 \begin{split}
|\zeta (Z')| \leq &\int_{Z'} J_\tau \zeta \, \d \Ha^n \\
&= \int_{\Sigma_1'} \int_0^{\tilde R} \prod_{i=1}^n |1 - tk_i(x)|  \, \d t \d \Ha^n +  \int_{\Sigma_2'} \int_0^{1/k_n(x)} \prod_{i=1}^n |1 - tk_i(x)|  \, \d t \d \Ha^n.
 \end{split}
\eeq
By the definition of $\Sigma_1'$ it holds $|1 - tk_i(x)| = (1 - tk_i(x))$ for every $(x,t)  \in \Sigma_1' \times [0,\tilde R]$ and therefore by the arithmetic-geometric inequality we may estimate 
\[
\prod_{i=1}^n |1 - tk_i(x)| \leq C(1 + |H_E(x)|^n) \qquad \text{for } \, (x,t)  \in \Sigma_1' \times [0,\tilde R].
\]
Similarly, we deduce that
\[
\prod_{i=1}^n |1 - tk_i(x)| \leq C(1 + t^n |H_E(x)|^n) \qquad \text{for } \, x   \in \Sigma_2' \ \ \text{and} \  0 \leq t \leq 1/k_n(x). 
\]
On the other hand, by the definition of $\Sigma_2'$ it holds $1/k_n(x) < \tilde R$. Therefore, by   \eqref{monros3},  $\lambda \leq C$  and \eqref{monros1} we have 
\[
 \begin{split}
|\zeta (Z')| &\leq C \int_{\Sigma_1' \cup \Sigma_2'} \int_0^{\tilde R} (1 + |H_E(x)|^n)    \, \d t \d \Ha^n  = C \tilde R \int_{\pa E \setminus \Sigma}  (1 + |H_E(x)|^n)    \, \d \Ha^n \\
&\leq C \int_{\pa E \setminus \Sigma}  (1+ \lambda^n  + |H_E - \lambda |^n)    \, \d \Ha^n  \leq C (\Ha^n(\pa E \setminus \Sigma) +  \| H_E - \lambda\|_{L^n(\pa E)}^n) \\
&\leq C \| H_E - \lambda\|_{L^n(\pa E)},
 \end{split}
\]
when $\| H_E - \lambda\|_{L^n(\pa E)} \leq 1$.  Hence by decreasing $\delta$, if needed, we have \eqref{monros3}. In particular, it holds
\beq
\label{monros4}
R_1  + R_2 + R_3 + R_4 \leq C \| H_E - \lambda\|_{L^n(\pa E)}
\eeq
where $R_i$ are  defined in \eqref{R1}-\eqref{R4}. 
\ \\
\ \\
\textbf{Step 2:} \quad Here we utilize the estimate \eqref{monros4} and prove the following auxiliary result. For a Borel set $\Gamma \subset \pa E$ and $0<r< R$ it holds
\beq
\label{monros5}
|E \cap \zeta (Z \cap (\Gamma \times (r,R)))| \geq \frac{ \Ha^n(\Gamma)}{(n+1) R^{n}} (R-r )^{n+1} - C\| H_E - \lambda\|_{L^n(\pa E)}.
\eeq

We prove \eqref{monros5} by 'backtracking' the Montiel-Ros argument. By the definition of  $R_1, R_2, R_3, R_4$  and \eqref{monros4} we may estimate  
\[
\begin{split}
|E \cap \zeta (Z \cap (\Gamma \times (r,R)))|  &\geq |\zeta (Z \cap (\Gamma \times (r,R)))| - R_1 \\
&\geq \int_{\zeta (Z \cap (\Gamma \times (r,R)))} \Ha^0(\zeta^{-1}(y) \cap Z\cap  (\Gamma \times (r,R))) \, dy  - R_1 -R_2 \\
&=  \int_{\Gamma \cap \Sigma}\int_{\min\{r, 1/k_n(x)\}}^{\min\{R, 1/k_n(x)\}} \prod_{i=1}^n (1 - tk_i(x))  \, \d t \d \Ha^n  - R_1 -R_2   \\
&\geq   \int_{\Gamma \cap \Sigma} \int_{\min\{r, 1/k_n(x)\}}^{\min\{R, 1/k_n(x)\}}  \left(1 - \frac{t}{n} H_E(x)\right)^n  \, \d t \d \Ha^n   - R_1 -R_2 - R_3 \\
&\geq  \int_{\Gamma \cap \Sigma} \int_{\min\{r, 1/k_n(x)\}}^{\min\{R, n/H_E(x)\}}\left(1 - \frac{t}{n} H_E\right)^n  \, \d t \d \Ha^n   - R_1 -R_2 - R_3-R_4 \\
&\geq  \int_{\Gamma \cap \Sigma} \int_{\min\{r, n/H_E(x)\}}^{\min\{R, n/H_E(x)\}}\left(1 - \frac{t}{n} H_E\right)^n  \, \d t \d \Ha^n   - R_1 -R_2 - R_3-R_4.
\end{split}
\]
Recall that for $x \in \Sigma$ it holds  $ \lambda /2 \leq H_E(x) \leq 2 \lambda$ and that $R = n/\lambda$. Therefore, we may estimate 
\[
\begin{split}
\int_{\Gamma \cap \Sigma} \int_{\min\{r, n/H_E(x)\}}^{\min\{R, n/H_E(x)\}}\left(1 - \frac{t}{n} H_E\right)^n  \, \d t \d \Ha^n  
&\geq \int_{\Gamma \cap \Sigma} \int_{\min\{r, n/H_E(x)\}}^{\min\{R, n/H_E(x)\}}\left(1 - \frac{t}{n} \lambda \right)^n  \, \d t \d \Ha^n  - C  \| H_E - \lambda\|_{L^n(\pa E)}  \\
&\geq   \int_{\Gamma \cap \Sigma} \int_{r}^{R}\left(1 - \frac{t}{n} \lambda  \right)^{n} \, \d t \d \Ha^n  - C  \| H_E - \lambda\|_{L^n(\pa E)} \\
&=   \frac{ \Ha^n(\Gamma \cap \Sigma) \, n}{(n+1) \lambda} \left(1 - \frac{\lambda}{n} r\right)^{n+1} - C  \| H_E - \lambda\|_{L^n(\pa E)} \\
&=   \frac{ \Ha^n(\Gamma \cap \Sigma) \, R}{(n+1)} \left(1 - \frac{r}{R}\right)^{n+1} - C  \| H_E - \lambda\|_{L^n(\pa E)}. 
\end{split}
\]
Hence, we obtain \eqref{monros5} from the previous two inequalities, from \eqref{monros1} and \eqref{monros4}. 
\newline
\newline
\textbf{Step 3:} \quad Here we finally prove the proposition. Recall the definition of  $E_r$  in \eqref{E_r}. Let us first prove that 
\beq
\label{monros6}
|E_r| \geq  \frac{P(E)}{(n+1) R^{n}} (R-r )^{n+1} - C\| H_E - \lambda\|_{L^n(\pa E)}
\eeq
for all $0<r<R$. 

To this aim, we claim that it holds
\beq
\label{monros7}
E \cap \zeta(Z \cap ( \Sigma \times (r,R))) \subset E_r  \cup   \{ y \in \zeta(Z) : \Ha^0(\zeta^{-1}(y) \cap Z) \geq 2 \}  \cup \zeta(Z').
\eeq 
The point of  this inclusion is  that    almost every point  which is of the form  $y = x -t \nu_E(x)$,  for $x \in Z$ and $t \in (r,R)$, belongs to $E_r$.

To this aim let $y \in E \cap \zeta(\Sigma \times (r,R))$. Then we may write  $y = x - t \nu_E(x)=\zeta(x,t)$ for some $x \in \Sigma$ and $t \in (r,R)$, with $(x,t) \in Z$. If $d_{\pa E}(y)= |y-x|$ then $y \in E_r$ because $|x-y| = t >r$. 
Otherwise, $d_{\pa E}(y) = |y-\tilde x| = \tilde r < t$ for $ \tilde x \in \pa E$, so we may write $y =  \tilde x - \tilde r \nu_E(x) =\zeta(\tilde x,\tilde r)$
 and $(\tilde x, \tilde r) \in Z \cup Z'$. Again, if $(\tilde x, \tilde r) \notin Z'$, then  $(\tilde x, \tilde r) \in Z$
and thus $\Ha^0(\zeta^{-1}(y) \cap Z) \geq 2$. Hence, we have \eqref{monros7}.

Recall that by the definition of $R_2$  and by  \eqref{monros4} it holds
\beq \label{monros71}
\begin{split}
|  \{ y \in \zeta(Z) :  \Ha^0(\zeta^{-1}(y) \cap Z) \geq 2 \}| & \leq \int_{\zeta(Z)} | \Ha^0(\zeta^{-1}(y) \cap Z) -1| \, \d y\\
&\leq C  \| H_E - \lambda\|_{L^n(\pa E)}.
\end{split}
\eeq
We then use \eqref{monros7}, \eqref{monros71},   \eqref{Z'} and   \eqref{monros5} with $\Gamma = \Sigma$  to deduce 
\[
\begin{split}
|E_r| &\geq |E \cap \zeta(Z \cap (\Sigma \times (r,R)))| -  C  \| H_E - \lambda\|_{L^n(\pa E)} \\
&\geq  \frac{\Ha^n(\Sigma)}{(n+1) R^{n}} (R-r )^{n+1} - C\| H_E - \lambda\|_{L^n(\pa E)}.
\end{split}
\]
The inequality \eqref{monros6} then follows from \eqref{monros1}.

Let us next  show that  for all $r \in (0,R)$ it holds 
\beq
\label{monros8}
|E_r| \leq \frac{\Ha^n(\Gamma_r)}{(n+1) R^{n}} (R-r)^{n+1} + C\| H_E - \lambda\|_{L^n(\pa E)},
\eeq
where $\Gamma_r \subset \pa E$ is defined in \eqref{gamma_r}.

First, we show
\beq
\label{monros9}
|E_R| \leq C\| H_E - \lambda\|_{L^n(\pa E)}.
\eeq
This follows from  an  already familiar argument, so we only sketch it. It is easy to see that   $E_R \subset \zeta(Z') \cup \zeta(Z \cap (\Sigma \times (R,\infty)))$.  Moreover, since $\lambda/2 \leq H_E(x) \leq 2 \lambda$ for $x \in \Sigma$, it holds 
\[
J_\tau \zeta(x,t) = \prod_{i=1}^n |1 - tk_i(x)| \leq C(1 + |H_E(x)|^n) \leq C \qquad \text{for } \, (x,t)  \in Z \cap (\Sigma \times (R,\infty)).  
\]
Recall that $R= n /\lambda$. Therefore, we have 
\[
\begin{split}
|\zeta(Z \cap (\Sigma \times (R,\infty))| &\leq \int_\Sigma \int_R^{\max\{n/H_E(x), R\}} J_\tau \zeta(x,t) \, \d t \Ha^n \\
&\leq C \int_\Sigma \big{|} \frac{n}{H_E} - R\big{|}\, \d t \Ha^n  \leq C \| H_E - \lambda\|_{L^n(\pa E)}. 
\end{split}
\]
The estimate \eqref{monros9} then follows from $|E_R| \leq |\zeta(Z \cap (\Sigma \times (R,\infty))| + |\zeta(Z')|$ and   \eqref{Z'}.

 Note that for all $\rho \in (r, R)$ it holds $\{ x \in E : d_{\pa E}(x) = \rho\} = \zeta(\Gamma_\rho, \rho)$  and $\Gamma_\rho \subset \Gamma_r$. We set $\zeta_\rho = \zeta(\cdot ,\rho): \pa E \to \R^{n+1}$ and thus it holds $\{ x \in E : d_{\pa E}(x) = \rho\} = \zeta_\rho(\Gamma_\rho)$  and
\[
J_\tau \zeta_\rho (x) = \prod_{i=1}^n |1 - \rho k_i(x)| \leq \left( 1- \frac{H_E}{n} \rho \right)^n \qquad \text{for } \, x  \in \Gamma_\rho.
\]
Therefore, by \eqref{monros9} and by co-area and area formulas we obtain
\[
\begin{split}
|E_r| &\leq  |E_r| - |E_R| + C \| H_E - \lambda\|_{L^n(\pa E)} \leq  \int_{r}^R \Ha^n(\{ x \in E : d_{\pa E} = \rho\} )\, \d\rho + C \| H_E - \lambda\|_{L^n(\pa E)}  \\
&=  \int_{r}^R \Ha^n(\zeta_\rho(\Gamma_\rho))\, \d\rho + C \| H_E - \lambda\|_{L^n(\pa E)}\\
&\leq \int_{r}^R \int_{\Gamma_\rho} J_\tau \zeta_\rho(x)\,\d\Ha^n  \d\rho+ C \| H_E - \lambda\|_{L^n(\pa E)}\\
&\leq \int_{r}^R \int_{\Gamma_\rho} \left( 1- \frac{H_E}{n} \rho \right)^n \,\d\Ha^n  \d\rho + C \| H_E - \lambda\|_{L^n(\pa E)}\\
&\leq  \int_{r}^R \Ha^n(\Gamma_\rho) \left( 1- \frac{\lambda}{n} \rho \right)^n\, \d\rho + C \| H_E - \lambda\|_{L^n(\pa E)}\\
&\leq  \Ha^n(\Gamma_r) \int_{r}^R  \left( 1- \frac{\rho}{R}  \right)^n\,  \d\rho + C \| H_E - \lambda\|_{L^n(\pa E)}\\
&= \frac{\Ha^n(\Gamma_r)}{(n+1) R^{n}} (R-r)^{n+1} + C\| H_E - \lambda\|_{L^n(\pa E)}.
\end{split}
\]
Hence, we have \eqref{monros8}.

 The second claim of the proposition follows immediately from \eqref{monros6} and \eqref{monros8}. These also imply 
\[
\Big{|}|E_r| -  \frac{P(E)}{(n+1) R^{n}} (R-r )^{n+1}\Big{|}  \leq  C\| H_E - \lambda\|_{L^n(\pa E)}.
\]
The first claim thus follows from \eqref{peri-vol} and $R = n/\lambda$.

For the last claim we refine the inclusion \eqref{monros7} and show that  for $0 < \rho < r <R$ and $r' \in (r,R)$   it holds
\beq
\label{monros10}
E \cap \zeta(Z \cap (\Gamma_{r'}  \times (r'-\rho,R))) \subset ( E_r + B_\rho)  \cup    \{ y \in \zeta(Z) : \Ha^0(\zeta^{-1}(y) \cap Z) \geq 2 \} \cup \zeta(Z').
\eeq 
Indeed, let $y \in E \cap \zeta(Z \cap (\Gamma_{r'}  \times (r'-\rho,R))) ) $.  Then we may write $y = x - t \nu_E(x)$ for some $x \in \Sigma \cap \Gamma_{r'}$ and $t \in ({r'}-\rho,R)$, with $(x,t) \in Z$.   If $t \in (r',R)$ then   by  \eqref{monros7} it holds 
\[
\begin{split}
y \in E \cap \zeta(Z \cap(\Sigma \times (r,R))) &\subset  E_{r} \cup   \{ y \in \zeta(Z) : \Ha^0(\zeta^{-1}(y) \cap Z) \geq 2 \}  \cup  \zeta(Z') \\
&\subset ( E_r + B_\rho)  \cup    \{ y \in \zeta(Z) : \Ha^0(\zeta^{-1}(y) \cap Z) \geq 2 \}  \cup  \zeta(Z').
\end{split}
\] 
Let us then assume that $t \in (r'-\rho,r']$. We  write $y =x - r' \nu_E(x) + (r'-t)\nu_E(x)$. Since $x \in \Gamma_{r'}$, i.e., $\pa E$ satisfies the interior ball condition at $x$ with radius $r'>r$, then necessarily $x - r'\nu_E(x) \in E_r$. Therefore, 
since $0\leq r'-t < \rho$, we conclude that $y \in E_r + B_\rho$ and \eqref{monros10} follows.

We use \eqref{Z'}, \eqref{monros5},  \eqref{monros71} and \eqref{monros10}   to conclude
\[
\begin{split}
|E_r+ B_\rho|  &\geq |E \cap \zeta(Z \cap(\Gamma_{r'} \cap \times ({r'}-\rho,R)))| -  C  \| H_E - \lambda\|_{L^n(\pa E)} \\
&\geq  \frac{ \Ha^n(\Gamma_{r'})}{(n+1) R^{n}} (R-({r'}-\rho) )^{n+1} - C\| H_E - \lambda\|_{L^n(\pa E)}.
\end{split}
\]
By using the second claim of the proposition and then letting $r' \to r$ we deduce 
\[
| E_r + B_\rho|  \geq \frac{P(E)}{(n+1) R^{n}}   (R-(r-\rho) )^{n+1} - \frac{C}{(R-r)^{n+1}} \| H_E - \lambda\|_{L^n(\pa E)}.
\]
 On the other hand, it clearly holds $E_r + B_\rho \subset E_{r-\rho}$. Then by \eqref{monros8}  we have 
\[
|E_r + B_\rho|  \leq |E_{r-\rho}| \leq \frac{P(E)}{(n+1) R^{n}} (R-(r-\rho))^{n+1} + C\| H_E - \lambda\|_{L^n(\pa E)}.
\]
The last claim thus follows from the two previous inequalities and \eqref{peri-vol}. 
\end{proof}

\subsection{Proof of Theorem \ref{thm1}}

\begin{proof}[Proof of Theorem \ref{thm1}]
Let $E$, $\lambda$, $C_0$ be as in the formulation of Theorem \ref{thm1}.
Recall that we denote $R=n/\lambda$. As before 
$C$ denotes a constant which may change from line to line but always depends only on $C_0$ and $n$. 
Let us denote 
\[
\eps := \| H_E - \lambda\|_{L^n(\pa E)}.
\]
 If  $\eps =0$, then $E$ is a disjoint union  balls by \cite{DM}. Let us then assume that $0<\eps \leq \delta $, where $\delta$ is initially set as in Proposition \ref{monros}. 
We might shrink \(\delta\) several times but always in such a way that it depends only on \(C_0\) and the dimension \(n\). 
Indeed, by shrinking $\delta$, if needed, Proposition \ref{lambda-diam control} provides the estimates 
\[
1/C \leq \lambda, R \leq C 
\]
and hence the first claim of Theorem \ref{thm1} is clear. We will use these estimates repeatedly without further mention.

By Proposition \ref{lambda-diam control} the number of 
the connected components of $E$ and their diameters are bounded by $C$. Thus, by applying a similar argument as in the proof of Proposition \ref{monros} (to obtain \eqref{peri-vol})
on each component and then summing these estimates up we obtain
\beq
\label{peri-vol2}
\big{|}nP(E) - (n+1)\lambda \,  |E|  \big{|} \leq C\eps.
\eeq
By possibly shrinking $\delta$ we have $R-\delta^{\frac{1}{n+2}}\geq  R/2$ . Choose $r_0 = R- \eps^{\frac{1}{n+2}}$. Then the volume estimates given by Proposition \ref{monros} read as
\beq \label{thm11}
\Big{|}  |E_r|  - \frac{|E|}{R^{n+1}} (R-r )^{n+1}  \Big{|} \leq  C \eps
\eeq
for all $ 0 \leq  r< R$ and
\beq \label{thm12}
\Big{|}  |E_r + B_\rho|  - \frac{|E|}{R^{n+1}} (R-(r-\rho) )^{n+1}  \Big{|} \leq  C \eps^{\frac{1}{n+2}}
\eeq
for all $0 \leq \rho \leq r \leq r_0$.
We remark that by \eqref{thm11}  we have 
\[
|E_{r_0}| \geq \frac{|E|}{R^{n+1}} \eps^{\frac{n+1}{n+2}} - C\eps \geq \frac{1}{C}\eps^{\frac{n+1}{n+2}} -C\eps.
\]
Hence by decreasing $\delta$, if needed, we may assume that  $E_{r_0}$  is non-empty.
This implies that $E_{r'}$ is non-empty for $r'>r_0$, when $|r'-r_0|$ is small enough. Since for any $r'>r_0$  it is geometrically clear that
 $\Gamma_{r'} \subset \partial E_{r_0} + \overline B_{r_0}$, then by using Proposition \ref{monros} and $r_0 = R-\eps^\frac{1}{n+2}$ we have
\[
\Ha^{n}(\pa E \setminus (\overline E_{r_0}+\overline B_{r_0})) \leq 
\Ha^{n}(\pa E \setminus \Gamma_{r'}) \leq C\frac{\eps}{(r_0-r' + \eps^\frac{1}{n+2})^{n+1}}.
\]
Thus by letting $r' \rightarrow r_0$ the previous estimate yields
\beq 
\label{thm13}
\Ha^{n}(\pa E \setminus (\overline E_{r_0}+\overline B_{r_0})) \leq C \eps^{\frac{1}{n+2}}.
\eeq
As previously, we divide the proof into three steps.  
\ \\
\ \\
\textbf{Step 1:} \quad
Recall that    $r_0 = R- \eps^{\frac{1}{n+2}} \geq R/2$.   We prove that there is a positive constant $d_0=d_0(C_0,n) \leq R/4$ such that  if $x,y \in E_{r_0}$, then 
\beq \label{thm14}
\text{either } \qquad |x-y| <  \eps^{\frac{1}{2(n+2)}} \qquad \text{or} \qquad |x-y| \geq d_0.
\eeq

Let us fix $x,y \in E_{r_0}$. We denote $d := |x-y|$ and  the segment from $x$ to $y$ by $J_{xy} := \{tx + (1-t)y :  t \in [0,1] \}$. We may assume that $d$ is small, since otherwise the claim \eqref{thm14} is trivially true. To be more precise we
 assume
\beq 
\label{thm14b}
d \leq \min\left\{\frac{R}{4},1\right\}.
\eeq

Let us first show that 
\beq \label{thm15}
J_{xy} \subset E_{r_0 - R^{-1}d^2}.
\eeq
Note that $ r_0 - R^{-1}d^2>0$ by $r_0\geq R/2$ and \eqref{thm14b} and hence $E_{r_0 - R^{-1}d^2}$ is well-defined and non-empty.
Choose  $z \in \mathbb R^{n+1} \setminus E$ and  $z' \in J_{xy}$ such that 
\[
|z-z'| = \dist(\mathbb R^{n+1} \setminus E,J_{xy}).
\]
If $z' = x$ or $z' = y$, then it follows from $x,y \in E_{r_0}$ that $|z-z'|> r_0 $. If not, then from the fact that  $z'$ is the closest point on $J_{xy}$ to $z$, we deduce that the vector $x-z'$ is orthogonal to $z-z'$, i.e., $\langle x-z',z-z'\rangle = 0$. Note also that 
 $\min \{|x-z'|, |y-z'| \} \leq d/2$ and we  may thus assume that $|x-z'| \leq d/2$. Therefore,  we have by Pythagorean theorem 
\[
|x-z|^2 = |x-z'|^2 + |z-z'|^2 \leq \frac{d^2}{4} + |z-z'|^2.
\]
Since $|x-z|> r_0 $, the previous estimate gives us
\[
|z-z'|^2 > r_0^2 - \frac{d^2}{4}. 
\]
We deduce from $r_0 \geq R/2$ and \eqref{thm14b} that
\[
 \left(r_0^2 - \frac{d^2}{4}\right)^{1/2} \geq r_0 - \frac{d^2}{R}.
\]
The previous two estimates yield $|z-z'| >  r_0 - R^{-1}d^2$ and the claim  \eqref{thm15} follows due to the choice of \(z\) and \(z'\).

Again, we use  $r_0 \geq R/2$ and \eqref{thm14b} to observe 
\[
r_0 - (1+R^{-1})d^2 \geq r_0 - d -  R^{-1}d^2 \geq \frac{R}{2}-\frac{R}{4}-\frac{R}{16}>0.
\]
Thus $E_{r_0 - (1+R^{-1})d^2}$ is well-defined and non-empty.
Next, we deduce from \eqref{thm15} and $E_r + B_\rho \subset E_{r-\rho}$ that 
\beq \label{thm16}
J_{xy} + B_{d^2}  \subset E_{r_0 - R^{-1}d^2} + B_{d^2} \subset E_{r_0 - (1+R^{-1})d^2}.
\eeq
Since $J_{xy} + B_{d^2}$ contains the cylinder $J_{xy} \times B_{d^2}^n$, it is clear that 
\[
|J_{xy} + B_{d^2} | \geq \omega_n d^{1+2n}.
\]
On the other hand, \eqref{thm11} and $\eps \leq 1$ (we may assume $\delta \leq 1$) imply
\begin{align*}
|E_{r_0 - (1+R^{-1})d^2}| &\leq \frac{|E|}{R^{n+1}} \big(R- (r_0 - (1+R^{-1})d^2)  \big)^{n+1} + C\eps \\
 &=\frac{|E|}{R^{n+1}} \big(\eps^\frac{1}{n+2} + (1+R^{-1})d^2)  \big)^{n+1} + C\eps \\
&\leq \frac{|E|}{R^{n+1}} \big(\eps^\frac{1}{n+2} + (1+R^{-1})d^2)  \big)^{n+1} + C\eps^\frac{n+1}{n+2} \\
&\leq C d^{2(n+1)} + C \eps^{\frac{n+1}{n+2}}.
\end{align*}
Then \eqref{thm16} yields
\[
\omega_n d^{1+2n}  \leq C d^{2(n+1)} + C \eps^{\frac{n+1}{n+2}}.
\]
If $d \geq \eps^{\frac{1}{2(n+2)}}$, then 
\[
\omega_n d^{1+2n}  \leq C d^{2(n+1)}.
\]
This implies $d \geq c >0$ for  some $c=c(C_0,n)$. By  recalling \eqref{thm14b}  the claim \eqref{thm14} follows. 
\newline
\newline
\textbf{Step 2:} \quad 
By \eqref{thm14} and possibly replacing $\delta$ with $\min\{\delta,(d_0/8)^{2(n+2)}\}$ we may divide the set $E_{r_0}$ into $N$ many clusters $E_{r_0}^1, \dots, E_{r_0}^N$ such that  we fix a point $x_i \in E_{r_0}$ and define the corresponding cluster  $E_{r_0}^i$ as
\[
E_{r_0}^i = \{ x \in E_{r_0} : |x-x_i| \leq d_0/8 \}. 
\] 
By \eqref{thm14} it holds $ E_{r_0}^i \subset B_{\eps_0} (x_i)$, where $\eps_0 = \eps^{\frac{1}{2(n+2)}}$, and  $|x_i-x_j| \geq d_0$ for $i \neq j$. 
Therefore, we have for every $\rho>0$
\beq \label{thm17}
\bigcup_{i=1}^N B_\rho(x_i) \subset  E_{r_0} + B_\rho \subset \bigcup_{i=1}^N B_{\rho+\eps_0} (x_i).
\eeq
Since $r_0 \geq R/2> R/4 \geq d_0$ and $|x_i-x_j| \geq d_0$ for $i \neq j$, then
 the balls $B_{\rho}(x_1), \ldots,B_{\rho}(x_N) $ with  $\rho = d_0/4$  are disjoint and contained in $E$, which, in turn, implies there is an upper bound $N_0 = N_0(C_0,n) \in \mathbb N$
for the number of clusters \(N\).

Next we  improve the lower bound  $|x_i -x_j| \geq d_0$  and  prove that there is a positive constant \(C_1=C_1(C_0,n)\)
such that
\beq \label{thm18}
|x_i-x_j| \geq 2R - 2C_1\eps^{\frac{1}{(n+2)^2}} \qquad  \text{for all pairs } \,  i \neq j.
\eeq
As a byproduct we  prove the last statement of the theorem, i.e., we show 
\beq
\label{thm18b}
\Big{|} P(E) -  N (n+1) \omega_{n+1} R^{n}  \Big{|} \leq  C \eps^{\frac{1}{2(n+2)}}.
\eeq

Recall that   the balls $B_{d_0/4}(x_1), \dots, B_{d_0/4}(x_N)$  are disjoint. Therefore, using \(N \leq N_0\) and \eqref{thm17} with  $\rho = d_0/4$  we deduce  
\[
\Big{|} |E_{r_0} + B_{d_0/4}| - N \omega_{n+1} \left(\frac{d_0}{4} \right)^{n+1} \Big{|} \leq C \eps_0 = C \eps^{\frac{1}{2(n+2)}}.
\]
On the other hand, we have
\(
d_0/4 \leq R/16 < R/2 \leq r_0
\)
so we may use \eqref{thm12} to obtain
\[
\Big{|} |E_{r_0} + B_{d_0/4}| - \frac{|E|}{R^{n+1}} \left(\frac{d_0}{4}  + \eps^{\frac{1}{n+2}} \right)^{n+1}  \Big{|} \leq C \eps^{\frac{1}{n+2}}.
\]
These two estimates and $\eps \leq 1$ imply
\beq \label{thm19}
\Big{|} |E| -  N \omega_{n+1} R^{n+1}  \Big{|} \leq  C \eps^{\frac{1}{2(n+2)}}.
\eeq
Thus,  \eqref{peri-vol2}, $R = n/\lambda$  and \eqref{thm19} yield \eqref{thm18b}.

To obtain \eqref{thm18}, let us assume that there is $0<h<R/2$  such that $|x_i -x_j| < 2R - 2h$ for some $i \neq j$.
This implies that  the balls  $B_R (x_i)$ and $B_R (x_j)$ intersect each other such that a set enclosed by  a spherical cap of height $h$ is included in their intersection. As the volume enclosed  by the spherical  cap of height $h$ has a lower bound  $c_n R^{n+1} h^{\frac{n+2}{2}}$, with some dimensional constant $c_n$, then there is $c=c(C_0,n)$ such that
\[
\Big{|}B_R (x_i) \cap B_R (x_j)  \Big{|} \geq  c h^{\frac{n+2}{2}}.
\]
We use the previous estimate as well as \eqref{thm12}, \eqref{thm17}, \eqref{thm19}, $\eps \leq 1$ and \(N \leq N_0\) to estimate
\begin{align*}
N\omega_{n+1}R^{n+1}
&\leq |E| + C \eps_0 \\
&\leq |E_{r_0}+B_{r_0}| + C \eps_0+ C \eps^\frac{1}{n+2} \\
&\leq \left| \bigcup_{i=1}^N B_{R+\eps_0} (x_i)\right| +C\eps_0 + C \eps^\frac{1}{n+2} \\
&\leq \left| \bigcup_{i=1}^N B_R (x_i)\right| + N\omega_{n+1} ((R+\eps_0)^{n+1} - R^{n+1}) +  C\eps_0 + C \eps^\frac{1}{n+2} \\
&\leq  N\omega_{n+1}R^{n+1} - \Big{|}B_R (x_i) \cap B_R (x_j)  \Big{|}+ C\eps_0 + C \eps^\frac{1}{n+2} \\
&\leq  N\omega_{n+1}R^{n+1} -  c h^{\frac{n+2}{2}}+ C\eps_0 + C \eps^\frac{1}{n+2} \\
&=  N\omega_{n+1}R^{n+1} -  c h^{\frac{n+2}{2}}+ C\eps^\frac{1}{2(n+2)} + C \eps^\frac{1}{n+2} \\
&\leq  N\omega_{n+1}R^{n+1} -  c h^{\frac{n+2}{2}}+ C \eps^\frac{1}{2(n+2)}.
\end{align*}
Thus $h^{\frac{n+2}{2}} \leq C \eps^\frac{1}{2(n+2)}$ and  \eqref{thm18} follows.
\newline
\newline
\textbf{Step 3:}\quad
Let $C_1$ be as in \eqref{thm18}.
By decreasing $\delta$, if needed, we may assume  
\[
0< R - C_1\eps^\frac{1}{(n+2)^2}<R-\eps^\frac{1}{n+2} =  r_0.
\]
Then we have by \eqref{thm17} and  \eqref{thm18} that the balls $B_\rho(x_1), \dots, B_\rho(x_N)$, with $\rho = R - C_1\eps^{\frac{1}{(n+2)^2}}$, are disjoint and 
\beq
\label{thm119}
\bigcup_{i=1}^N B_\rho(x_i) \subset E_{r_0} + B_\rho  \subset E_{r_0 - \rho} \subset E.
\eeq
This, $\eps \leq 1$, $N \leq N_0$ and  \eqref{thm19} imply  
\beq \label{thm120}
\left | E \setminus  \bigcup_{i=1}^N B_\rho(x_i)  \right | \leq C \eps^{\frac{1}{(n+2)^2}}.
\eeq

Set $\eps_1= \eps^{\frac{1}{(n+2)^3}}$. We prove 
\beq \label{thm121}
E \subset \bigcup_{i=1}^N B_{\eta}(x_i)
\eeq
for $\eta= R + C_2 \eps_1$ with some positive $C_2=C_2(n,C_0)$. 
By decreasing $\delta$, if necessary, we deduce from \eqref{thm120} that
\[
|B_{\eps_1}| > \left | E \setminus  \bigcup_{i=1}^N B_\rho(x_i)  \right |.
\] 
Thus, if $x \in E_{\eps_1}$, then $B_{\eps_1}(x) \cap \bigcup_{i=1}^N B_\rho(x_i) $ must be non-empty. This implies
\beq \label{thm122}
E_{\eps_1} \subset  \bigcup_{i=1}^N B_{\rho+\eps_1}(x_i).
\eeq

Assume that  for $x \in \pa E$ it holds
\[
d_x : = \dist \left(x, \overline E_{r_0} + \overline B_{r_0}\right)>0.
\]
Then by \eqref{thm13}
\[
\Ha^n(\pa E \cap B(x,d_x)) \leq C \eps^\frac{1}{n+2}.
\]
Let $\delta_n \in \mathbb R_+$ be as in Lemma \ref{density} and set
$
r_x=\min\left\{d_x,\delta_n/\lambda\right\}.
$
Again, by possibly decreasing  $\delta$ so that $\delta \leq \delta_n$,  Lemma \ref{density} yields  
\[
\delta_n r_x^n \leq \Ha^n\left(\pa E \cap B_{r_x}(x)\right).
\]
By combining the two previous estimates we have
\[
\min\left\{d_x,\frac{\delta_n}{\lambda}\right\} \leq C\eps^\frac{1}{n(n+2)}.
\]
Since $\delta_n/\lambda \geq \delta_n/C$, then by decreasing $\delta$, if necessary, the previous estimate implies $r_x=d_x$ and further gives us
\beq
\label{thm123}
d_x \leq   C\eps^\frac{1}{n(n+2)} \leq C\eps^\frac{1}{(n+2)^2}.
\eeq
On the other hand, by \eqref{thm17}  
\beq
\label{thm124}
\overline E_{r_0} + \overline B_{r_0} \subset  E_{r_0} +  B_R \subset   \bigcup_{i=1}^N B_{R+\eps_0} (x_i),
\eeq
where $\eps_0=\eps^\frac{1}{2(n+2)} \leq \eps^\frac{1}{(n+2)^2}$. Thus, \eqref{thm123} and \eqref{thm124} imply
\[
\pa E \subset  \bigcup_{i=1}^N B_{\tilde \eta}(x_i)
\]
with $\tilde \eta = R+C\eps^\frac{1}{(n+2)^2}$. By combining this observation with \eqref{thm122} we obtain \eqref{thm121}.

Finally, by decreasing $\delta$ one more time, if necessary, \eqref{thm18b},  \eqref{thm119} and \eqref{thm121} yield 
\[
\bigcup_{i=1}^N B_{\rho_-} (x_i) \subset  E \subset \bigcup_{i=1}^N B_{\rho_+}(x_i), 
\]
where $\rho_-= R - C \eps^\frac{1}{(n+2)^3}$, $\rho_+ = R + C \eps^{\frac{1}{(n+2)^3}}$, the balls $B_{\rho_-} (x_1), \ldots,B_{\rho_-} (x_N)$
are mutually disjoint, for $N$ it holds
\[
\Big{|} P(E) -  N (n+1) \omega_{n+1} R^{n}  \Big{|} \leq  C \eps^{\frac{1}{(n+2)^3}}
\]
and $C=C(C_0,n) \in \R_+$. The claim of Theorem \ref{thm1} then follows by  Remark \ref{rem thm1}. 
\end{proof}


\section{Asymptotic behavior of the volume preserving mean curvature flow}

In this section we first define the flat flow and recall some of its  basic properties. We do this in the general dimensional case  $\R^{n+1}$ and resctrict ourself to the case $n \leq 2$ only in the proof of Theorem \ref{thm2}.  We  begin by  defining the flat flow of \eqref{the flow}.

Assume that $E_0 \subset \R^{n+1}$ is a bounded set of finite perimeter  with the volume of the unit ball   $|E_0| = \omega_{n+1}$. For given $h \in \R_+$ we construct a sequence of sets
$(E_k^h)_{k=1}^\infty$ by iterative minimizing procedure called minimizing movements, where initially $E_0^h = E_0$ and  $E_{k+1}^h$ is a minimizer of the following problem
\beq
\label{min mov}
\mathcal{F}_h(E, E_k) = P(E) + \frac{1}{h} \int_E \bar{d}_{E_k} \, dx  + \frac{1}{\sqrt{h}}\big| |E| - \omega_{n+1} \big|.
\eeq
Recall, that $\bar{d}_{E_k}$ is the signed distance function from $E_k$. We then define the approximative flat flow $(E_t^h)_{t \geq  0}$ by
\beq
\label{app flow}
E_t^h = E_k^h, \qquad \text{for } \,  (k-1)h \leq t < k h. 
\eeq
By \cite{MSS} we know that there is a subsequence of the approximative flat flow which converges
\[
(E_t^{h_l})_{t \geq  0} \to (E_t)_{t \geq 0},
\]
where for every $t >0$ the set $E_t$ is a set of finite perimeter with $|E_t| = \omega_{n+1}$. Any  such limit is called a flat flow of \eqref{the flow}.   It follows from the result in  \cite{CN} that  if the initial set $E_0$  is  smooth, then any flat flow coincide with the  unique classical solution of \eqref{the flow} as long as the latter is defined (see also \cite{CMP}).

\subsection{Preliminary results} 
Let us take more rigorous approach to the concepts heuristically introduced above. We base this mainly on \cite{MSS}, where
the only difference is that the volume constraint has a different value. Obviously, this does not affect the arguments. 

First, we take a closer look at the functional $\mathcal{F}_h$ given by \eqref{min mov}.
If $E,F \subset \mathbb R^{n+1}$ are bounded sets of finite perimeter, then it is easy to see that modifications of $E$ in a set of measure zero
do not affect the value $\mathcal{F}_h(E,F)$ whereas such modifications of $F$ may lead drastic changes of the of $\mathcal{F}_h(E,F)$.
To eliminate this issue, we use a convention that a topological boundary of a set of finite perimeter is always the support of the corresponding Gauss-Green measure. 
Thus, we consider $\mathcal{F}_h$ as a functional $X_{n+1} \times \{ A \in X_{n+1}: A \neq \varnothing\} \rightarrow \R$, where
\[
X_{n+1} = \{E \subset \R^{n+1}: E \ \ \text{is a bounded set of finite perimeter with} \ \ \pa E = \spt \ \mu_E\}.
\]
We remark that if $E_0$  is essentially open or closed and $E_0 \in X_{n+1}$, then we may assume it to be open or closed respectively. 

For a non-empty $F \in X_{n+1}$ there is always a minimizer $E$ of the functional $\mathcal F_h( \ \cdot \ , F)$ in the class $X_{n+1}$ satisfying the \emph{discrete dissipation inequality}
\beq
\label{dissipation}
P(E) + \frac{1}{h} \int_{E \Delta F} d_{\pa F} \ \d x + \frac{1}{\sqrt h} \left||E|-\omega_{n+1}\right| \leq P(F) +  \frac{1}{\sqrt h} \left||F|-\omega_{n+1}\right|,
\eeq
see \cite[Lemma 3.1]{MSS}.
Moreover, there is a dimensional constant $C_n$ such that 
\beq
\label{distance}
\sup_{E \Delta F} d_{\pa F} \leq C_n \sqrt{h},
\eeq
see \cite[Proposition 3.2]{MSS}. The minimizer $E$ is always a $(\Lambda,r_0)$ -\emph{minimizer} in any open neighborhood of $E$ with suitable $\Lambda, r_0 \in \R_+$ satisfying
$\Lambda r_0 \leq 1$. Thus, by the standard regularity theory \cite[Thm 26.5 and Thm 28.1]{Ma} $\pa^* E$ is relatively open in $\pa E$ and $C^{1,\alpha}$-regular with any 
$0<\alpha <1/2$ and the Hausdorff dimension of the singular part $\pa E \setminus \pa^* E$ is at most $n-7$. These imply, via the relative isoperimetric inequality,
that $E$ can always  be chosen as an open set. On the other hand, if $E$ is non-empty, it has a Lipschitz-continuous distributional mean curvature $H_E$ satisfying the Euler-Lagrange equation
\beq
\label{euler}
\frac{\bar d_F}{h} = -H_E + \lambda_E,
\eeq
where the \emph{Lagrange multiplier} can be written in the case $|E| \neq \omega_{n+1}$ as
\beq
\label{euler2}
\lambda_E = \frac{1}{\sqrt h} \mathrm{sgn} \left(\omega_{n+1}-|E|\right),
\eeq
see \cite[Lemma 3.7]{MSS}. Thus, by using standard elliptic estimates one can show that $\pa^*E$ is in fact $C^{2,\alpha}$-regular and \eqref{euler} holds
in the classical sense on $\pa^* E$. In particular, $E$ is $C^{2,\alpha}$-regular when $n \leq 6$. Moreover, if $x \in \pa E$ satisfies exterior or interior ball condition with any $r$, then it must belong to the reduced boundary  of $E$. This is well-known and follows essentially from  \cite[Lemma 3]{DM}.

Let us turn our focus back on flat flows. 
Let $E_0 \in X_{n+1}$  be a set with volume $\omega_{n+1}$ and $0<h<(\omega_{n+1}/P(E_0))^2$. Then we find a minimizer $E^h_1 \in X_{n+1}$ for $\mathcal F_h( \ \cdot \ , E_0)$
and by \eqref{dissipation} we have
$
\left||E^h_1|-\omega_{n+1}\right| \leq \sqrt h P(E_0) 
$ 
implying, via the condition $h<(\omega_n/P(E_0))^2$, that $E^h_1$ is non-empty. Again, we find
a minimizer $E^h_2 \in X_{n+1}$ for $\mathcal F_h( \ \cdot \ , E_1)$ and using \eqref{dissipation} twice we obtain 
$
\left||E^h_2|-\omega_{n+1}\right| \leq \sqrt h P(E_0) 
$ 
and thus $E^h_2$ is also non-empty. By continuing the procedure we find non-empty sets
$E^h_0,E^h_1,E^h_2,\ldots \in X_{n+1}$ as mentioned earlier, i.e., $E^h_0 = E_0$ and  $E^h_k$ is a minimizer of
$\mathcal F_h( \ \cdot \ , E_{k-1})$ for every $k \in \N$. Thus, we may define  an approximate flat flow $(E^h_t)_{t\geq 0}$, with the initial set $E_0$, defined by \eqref{app flow}.
Further a flat flow as a limit is defined as before. By iterating \eqref{dissipation}
we obtain
\beq
\label{iterated}
P(E_{kh}^h) + \frac{1}{h} \sum_{j=1}^k \int_{E_{jh}^h \Delta E_{(j-1)h}^h} d_{ \pa E_{(j-1)h}^h} \ \d x + \frac{1}{\sqrt h}||E_{kh}^h|-\omega_{n+1}| \leq P(E_0)
\ \ \text{for every $k \in \N$. }
\eeq
By the earlier discussion we may assume that $E^h_t$, for every $t \geq h$, is  an open set and $\pa E^h_t$ is  $C^2$-regular up to the singular part $\pa E^h_t \setminus \pa^* E^h_t$ with Hausdorff dimesion  at most $n-7$.
We use the shorthand notation $\lambda^h_t$ for the corresponding Lagrange multiplier.

Next, we list some basic properties of the approximative flat flow.
\begin{proposition}
\label{spadaro0}
Let $(E_t^h)_{t \geq 0}$ be an approximative flat flow starting from $E_0 \in X_{n+1}$ with volume $\omega_{n+1}$ and  $P(E_0) \leq C_0$. 
There is a positive constant $C=C(C_0,n)$ such that the following holds for every $0<h<(\omega_n/P(E_0))^2$:
\begin{itemize}
\item[(i)] For every $s,t$ with $h \leq s \leq t -h$  it holds $|E_s^h \Delta E_t^h| \leq C \sqrt{t-s}$. 
\item[(ii)] Suppose that for given $T_1 \geq 0$ it holds $|E_{T_1}^h|=\omega_{n+1}$. Then $P(E_{T_1}^h) \geq P(E_t^h)$ for every $t \geq T_1$ and
\[
\int_{T_1+h}^{T_2} \int_{\pa^* E_{t}^h} (H_{E_{t}^h} - \lambda_t^h)^2\, d \Ha^n \d t \leq C (P(E_{T_1}^h) - P(E_{T_2}^h))
\]
 for every $T_2 \geq T_1+h.$ Moreover, for every $h \leq T_1 < T_2$ it holds
\[
\int_{T_1}^{T_2} \int_{\pa^* E_{t}^h} (H_{E_{t}^h} - \lambda_t^h)^2\, d \Ha^n \d t \leq C P(E_0).
\]

\item[(iii)] For every $T>0$ there is $R=R(E_0,T)$ such that $E_t^h \subset B_R$ for all $0 \leq t \leq T$.
 
\item[(iv)] If $(h_k)_k$ is a sequence of positive numbers converging to zero, then up to a subsequence 
there exist  approximative flat flows $((E^{h_k}_t)_{t\geq 0})_k$ which converges to a flat flow $(E_t)_{t\geq 0}$, where $E_t \in X_{n+1}$,  in $L^1$-sense  in space and pointwise time, i.e., for every  $t\geq 0$ it holds 
\[
\lim_{h_k \to 0} |E_t^{h_k} \Delta E_t| = 0. 
\]  
The limit flow also satisfies $|E_s \Delta E_t| \leq C \sqrt{t-s}$ for every $0<s<t$ and $|E_t|= \omega_{n+1}$ for every $t \geq 0$.
\item[(v)] If $E_0$ is either open or closed, then the sequence in (iv) converges to $(E_t)_{t\geq 0}$ in $L^1$  in space and  compactly uniformly in time, i.e., for a fixed  $T$ it  holds 
\[
\lim_{h_k \to 0} \sup_{t \in [0,T]} |E_t^{h_k} \Delta E_t| = 0. 
\]  
Moreover, $|E_s \Delta E_t| \leq C \sqrt{t-s}$ for every $0\leq s<t$.
\end{itemize}
\end{proposition}

\begin{proof} The claims (i) - (iv) are essentially proved in \cite{MSS}, see the proofs of Proposition 3.5, Lemma 3.6 and Theorem 2.2.

To prove (v), we 
first show that
\[
|E^h_h \Delta E_0| \rightarrow 0 \ \ \text{as} \ \ h \rightarrow 0
\]
which immediately implies via (iv) that $|E_0 \Delta E_t| \leq C\sqrt t$ for every $t\geq0$ and hence the second claim of (v) holds.
Then the compactly uniform convergence in time is  a rather direct consequence of this and (i).

To this aim, let $(h_k)_k$ be an arbitrary sequence 
of positive numbers converging to zero. By (iii) and by the standard compactness property of sets of finite perimeter there is a bounded set of finite perimeter
$E_\infty$ such that, up to extracting a subsequence, $E_{h_k}^{h_k} \rightarrow E_\infty$ in $L^1$ -sense. In particular, by \eqref{iterated} we have
$|E_\infty|= \omega_{n+1} = |E_0|$. Again, by using $|E^{h_k}_{h_k} \Delta E_\infty| \rightarrow 0$ and \eqref{distance} we have
\[
|E_\infty \setminus \{y \in \mathbb R^n : \bar d_{E_0}(y) \leq  j^{-1}\}| = 0 \quad \text{and} \quad  |\{y \in \mathbb R^n : \bar d_{E_0}(y) \leq - j^{-1}\} \setminus E_\infty| = 0
\] 
for every $j \in \mathbb N$. Thus, by letting $j \rightarrow \infty$ we obtain $|E_\infty \setminus \bar {E_0}|=0$ and $|\mathrm{int}(E_0) \setminus E_\infty|=0$.
Since $E_0$ is open or closed, this means either $|E_\infty \setminus E_0|=0$ or $|E_0 \setminus E_\infty|=0$.
But now $|E_\infty| = |E_0|$ so the previous yields $|E_\infty \Delta E_0| =0$. Thus, $|E^{h_k}_{h_k} \setminus E_0| \rightarrow 0$ up to a subsequence 
and since $(h_k)_k$ was arbitrarily chosen, then it holds $|E^h_h \Delta E_0| \rightarrow 0$.  
\end{proof}

We note that the claim (v) does not hold for every bounded set of finite perimeter $E_0$. As an example one may construct  a wild set of finite perimeter $E_0$ such that $|E^h_h \Delta E_0| \geq c_0 >0$ for all $h>0$

By \cite[Corollary 3.10]{MSS} it holds, for a fixed time $T\geq h$, that the integral $\int_h^T |\lambda^h_t|^2 dt$ is uniformly bounded in $h$
and hence, via \eqref{euler2}, it holds  $\left|\{ t \in (h,T): |E^h_t|\neq \omega_{n+1}\}\right| \leq C h$, where $C$ depends also on $T$.
 We may improve this by using  Lemma \ref{lambda-diam control}.

\begin{proposition}
\label{lambda bound}
Let $C_0>0$ and $E_0 \in X_{n+1}$ be a set of finite perimeter with volume $\omega_{n+1}$ and $P(E_0) \leq C_0$. 
There are positive constants $C=C(C_0,n)$ and $h_0=h_0(C_0,n)$ such that if $h \leq h_0$ and $(E_t^h)_{t \geq 0}$ is an approximative flat flow starting from $E_0$, then for every $h\leq T_1  \leq T_2$ 
\begin{align*}
\int_{T_1}^{T_2} |\lambda^h_t|^2 dt &\leq C (T_2 - T_1 +1) \ \ \text{and} \\
\left|\{ t \in (T_1,T_2): |E^h_t|\neq \omega_{n+1}\}\right| & \leq Ch (T_2 - T_1 +1).
\end{align*}
\end{proposition}

\begin{proof}
 By \eqref{iterated} we may choose $h_0=h_0(C_0,n)$  such that $|E^h_t| \geq \frac{\omega_{n+1}}{2}$ whenever $h \leq h_0$.
We may also assume $C_0 > 2 \omega_{n+1}$ so $|E^h_t| \geq 1/C_0$ for $h \leq h_0$. Thus, by Lemma \ref{lambda-diam control} and $P(E^h_t) \leq C_0$ we find a positive $C=C(C_0,n)$ such that for every $t \geq h$ and $h \leq h_0$ it holds
\[
|\lambda^h_t|^2 \leq C\left(1+\int_{\pa^* E_{t}^h} (H_{E_{t}^h} - \lambda_t^h)^2 \d \Ha^n \right).
\]
Therefore, 
\[
\int_{T_1}^{T_2} |\lambda^h_t|^2 dt \leq C (T_2 - T_1) + C \int_{T_1}^{T_2} \int_{\pa^* E_{t}^h} (H_{E_{t}^h} - \lambda_t^h)^2 \d \Ha^n \d t.
\]
By   Lemma \ref{spadaro0} (ii)     we obtain the first inequality.
The first inequality implies, via \eqref{euler2}, the second inequality with the same constant $C$.
\end{proof}

We need also the following  comparison result for the proof. 
\begin{lemma}
\label{comparison}
Let $1 \leq C_0<\infty$. Assume $E_0 \in X_{n+1}$ is a set of finite perimeter with volume $\omega_{n+1}$ and $P(E_0) \leq C_0$, 
and  let  $F = \bigcup_{i=1}^N B_r(x_i)$ with $|x_i - x_j | \geq 2 r$ and $1/C_0 \leq r \leq C_0$. There is a positive constant
$\eps_0=\eps_0(C_0,n)$  such that if
$(E_t^h)_{t \geq 0 }$ is an approximative flat flow starting from $E_0$ and
\[
\sup_{x \in E_{t_0}^h \Delta F} d_{\pa F}(x) \leq \eps \ \ \text{with} \ \ \eps \leq \eps_0
\]
for $t_0 \geq 0$, then it holds 
\[
\sup_{x \in E_t^h \Delta F} d_{\pa F}(x) \leq C \eps^{\frac19} \qquad \text{for all } \, t_0 < t < t_0+\sqrt{\eps}
\]
provided that $h \leq \min\{\sqrt \eps, h_0\}$, where  $h_0=h_0(C_0,n)$ is as in Proposition \ref{lambda bound}.
\end{lemma}

\begin{proof}
Our standing assumptions are  $h \leq \min\{\sqrt \eps, h_0\}$ and  $\eps \leq \min \{1/(2C_0),1\}$. As usual, $C$ denotes
a positive constant which may change from line to line but depends only on the parameters $C_0$ and $n$.

Without loss of generality we may assume $t_0=0$.
Fix an arbitrary $x_i \in \{ x_1, \dots, x_N\}$. Up to translating the coordinates we may  assume that $x_i = 0$. 
We set for every $k=0,1,2,\ldots$
\[
\rho_k = \inf \{ |x| : x \in \R^{n+1} \setminus  E_{kh}^h\}  \quad  \text{and} \quad  r_k  = \min \{r,\rho_0,\ldots,\rho_k\}. 
\]
We claim that  it holds
\beq
\label{discrete est}
r_{k+1}^2 - r_k^2 \geq - C_1(1 + |\lambda^h_{(k+1)h}|) h,
\eeq
with some positive constant $C_1=C_1(C_0,n).$
 First,  if $r_{k+1} =  r_k$, the claim \eqref{discrete est} is trivially true. Thus, we may assume $r_{k+1} <  r_k$ which implies 
$\rho_{k+1} = r_{k+1} < r_k \leq \rho_k$. Then $\rho_k>0$ which, in turn, means
\[
\rho_k = \min_{\pa E^h_{kh}} |x|.
\]
Since $E_{(k+1)h}^h$ is bounded and open, there is a point $x \in \R^{n+1} \setminus  E_{(k+1)h}^h $ with $\rho_{k+1}= |x|$.
 Let $x'$ be a closest point to $x$ on $\pa E_{kh}^h$. Then
\[
 r_{k+1} +|\bar{d}_{E_{kh}^h}(x)|= |x| + |\bar{d}_{E_{kh}^h}(x)| \geq |x'| \geq \rho_k \geq r_k.
\]
The condition $|x|<\rho_k$ means $x \in E_{kh}^h$ so the previous estimate yields
\beq
\label{auxest1}
r_{k+1}-r_k \geq  \bar d_{ E_{kh}^h}(x).
\eeq
Again,
 $x \in E_{kh}^h\setminus E_{(k+1)h}^h$ so  by Lemma \ref{distance} $|\bar{d}_{E_{kh}^h}(x)| \leq C_n \sqrt h$
and hence
\beq
\label{auxest2}
r_{k+1}-r_k \geq -C_n\sqrt h.
\eeq

We split the argument  into two cases. First, if $r_{k+1} < C_n \sqrt h$, then by \eqref{auxest2} we have $r_k < 2C_n \sqrt h$.
Therefore, using \eqref{auxest2} we obtain 
\beq
\label{auxest3}
r_{k+1}^2-r_k^2 \geq  -C_n(r_{k+1}+r_k )\sqrt h \geq -3C_n^2 h.
\eeq
If $r_{k+1} \geq C_n \sqrt h$, then by \eqref{auxest2} $r_k \leq 2r_{k+1}$. 
Since $r_{k+1}>0$, then it holds  $x \in \pa E_{(k+1)h}^h$ and $E_{(k+1)h}^h$ satisfies interior ball condition of radius $r_{k+1}$ at $x$. 
Thus, by discussion in Section 2  $x$ belongs to the reduced boundary of $E_{(k+1)h}^h$  and  therefore by the maximum principle it holds $H_{ E_{(k+1)h}^h}(x) \leq  \frac{n}{r_{k+1}}$. Again, by the previous estimate, \eqref{auxest1}, the Euler-Lagrange equation  \eqref{euler} and $ r_{k+1} \leq C_0$  we obtain
\[
\frac{r_{k+1} - r_k}{h} \geq \frac{\bar d_{ E_{kh}^h}(x)}{h}  \geq  -\frac{n}{r_{k+1}} - |\lambda^h_{(k+1)h}| \geq- \frac{1}{r_{k+1}} \left( n + C_0|\lambda^h_{(k+1)h}| \right) .
\]
Therefore,
\beq
\label{auxest5}
\frac{r_{k+1}^2 - r_k^2}{h} 
\geq  - \left(1+\frac{r_k}{r_{k+1}}\right) \left( n + C_0|\lambda^h_{(k+1)h}| \right) 
\geq  - 3\left( n + C_0|\lambda^h_{(k+1)h}| \right).
\eeq
Thus, \eqref{auxest3} and \eqref{auxest5} yield the claim \eqref{discrete est} in the case $r_{k+1} < r_k$.

We iterate \eqref{discrete est} up to $K \in \N$, which is chosen so that $Kh \in (\sqrt{\eps}, 2\sqrt{\eps})$ (recall $h < \sqrt \eps$),
and use Proposition \ref{lambda bound} to obtain 
\beq
 \label{auxest6}
\begin{split}
r_K^2 - r_0^2 
&\geq - C_1 \sum_{k=0}^{K-1} (1+|\lambda^h_{(k+1)h}|)h \\
&=-C_1Kh -C_1 \int_h^{(K+1)h} |\lambda^h_t| \, dt \\
&\geq -2C_1\sqrt\eps -C_1 \int_h^{3\sqrt \eps} |\lambda^h_t| \, dt \\
&\geq -2C_1\sqrt\eps - \int_h^{3\sqrt \eps} \eps^{-\frac14} + \eps^{\frac14} |\lambda^h_{t}|^2 \, dt \\
&\geq -C \eps^{\frac14} \left(1+\int_h^{3\sqrt \eps}|\lambda^h_{t}|^2 \, dt\right) \geq -C \eps^{\frac14}.
\end{split}
\eeq
By the assumption $\sup_{x \in E_0 \Delta F} d_{\pa F}(x) \leq \eps$ we have  $r-\eps \leq r_0$.
Thus we divide $r_K^2 - r_0^2$  by $r_K + r_0$ and use $r_0 \geq r-\eps \geq r/2 \geq 1/(2C_0)$ as well as \eqref{auxest6} to find a positive constant $C_2=C_2(C_0,n)$  
such that $r_K \geq r - C_2\eps^{\frac14}$. This means that 
\[
\inf_{\R^{n+1}  \setminus E_{t}^h} \bar d_{B_r(x_i)}  \geq - C_2 \eps^{\frac14}\qquad \text{for all } \, t < \sqrt{\eps}
\]
and again due to the arbitrariness of $x_i \in \{x_1,\ldots,x_N\}$
\[
\inf_{ \R^{n+1}  \setminus E_{t}^h} \bar d_F \geq -C_2 \eps^{\frac14}\qquad \text{for all } \, t < \sqrt{\eps}. 
\]

To conclude the proof, we show that there is a positive constant $\eps_1=\eps_1(C_0,n)$ such that
\beq \label{compa 2}
 \sup_{E_{t}^h} \bar d_F \leq 2 \eps^{\frac19}\qquad \text{for all } \, t < \sqrt{\eps}
\eeq
provided that $\eps \leq \eps_1$.
To this aim we choose an arbitrary $x_0 \in \R^{n+1}  \setminus \bar F$ with $\bar d_F(x_0) \geq 2\eps^{\frac19}$.
We set for every $k=0,1,2,\ldots$
\[
\rho_k = \inf_{x \in  E_{kh}^h} |x-x_0| \ \ \text{and} \ \ r_k  = \min \{2\eps^{\frac19},\rho_1,\ldots,\rho_k\}.
\]
In particular, $r_k \leq 2C_0^\frac{1}{9}$. A slight modification of the procedure we used to obtain  \eqref{auxest6} yields
\[
r_K^2 - r_0^2  \geq - C \eps^{\frac14},
\]
where $K$ is the same as earlier. Again, the conditions $\sup_{x \in E_0 \Delta F} d_{\pa F}(x) \leq \eps$ and $\eps \leq 1$ imply  $ r_0 \geq 2\eps^{\frac19}-\eps \geq \eps^{\frac19}$. 
Thus
\[
r_K-r_0 \geq -C \frac{\eps^{\frac14}}{r_0} \geq -C \eps^\frac{5}{36} = - C \eps^\frac{1}{36} \eps^{\frac19}
\]
and thus $r_K \geq (1-C\eps^\frac{1}{36}) \eps^{\frac19}>\frac12  \eps^{\frac19}$, when $\eps$ is small enough.  Since  $x_0$, with $d_F(x_0) \geq 2\eps^{\frac19}$, was arbitrarily chosen we deduce that 
\[
E_{kh}^h \subset \{ x \in \R^{n+1} : d_F(x) \leq 2\eps^{\frac19} \} \qquad \text{for all } \, k =0, \dots, K. 
\]
 The claim \eqref{compa 2} then  follows from  the choice of $K$.

\end{proof}

\subsection{Proof of Theorem \ref{thm2}}

The proof of Theorem \ref{thm2}  is based on Theorem \ref{thm1}. We first use it together with the dissipation inequality in Proposition \ref{spadaro0} (ii) to deduce that there exists a sequence of times $t_j \to \infty$ such that the sets $E_{t_j}$ are close to a disjoint union of balls. Since perimeter of the approximative flat flow is essentially  decreasing  then the number of the balls is also monotone. In particular, we deduce that after some time, the sets $E_{t_j}$ are close to a fixed number, say $N$,  of balls. We use the second statement of Theorem  \ref{thm1} to deduce  that the perimeters of $E_{t_j}$ converges to the perimeter  of $N$ many balls with volume $\omega_{n+1}$ and thus the right-hand-side of the dissipation inequality   converges to zero. This allows us to  improve  our estimate and use Theorem \ref{thm1} again to deduce that the flat flow $E_t$ is close to a disjoint union of $N$ many balls for all  large $t$ except a set of times  with small measure. The statement then finally follows from Lemma \ref{comparison}.

\begin{proof}[Proof of Theorem \ref{thm2}]
 Assume that the initial set $E_0 \in X_{n+1}$ has the volume of the unit ball $|E_0|=\omega_{n+1}$, fix a positive 
$C_0$ with $C_0 \geq \max\{1,P(E_0)\}$ and assume $h < (C_0/\omega_{n+1})^2$. Let $(E_t)_{t\geq 0}$ be a flat flow starting from $E_0$ and let $(E_t^{h_l})_{t\geq 0}$ be an approximative flat flow which by Proposition \ref{spadaro0} converges to $(E_t)_{t\geq 0}$ locally uniformly in $L^1$.  We simplify the notation and denote the converging subsequence again by $h$.   Since we are now in the dimensions 2 and 3 ($n=1,2$), the sets $E^h_t$ are  $C^2$-regular.

\medskip

\textbf{Step 1:}\quad Let us denote 
\beq \label{def sigma}
\Sigma_h := \{ t \in (0,\infty) :  |E^h_t|\neq \omega_{n+1} \} 
\eeq
By  \eqref{iterated} and Proposition \ref{lambda bound}   we find a constant $h_0=h_0(C_0,n) < 1$ such that
$|E^h_t| \geq 1/C_0$ for every $t\geq 0$ and
\[
\left| (T_1, T_2) \cap \Sigma_h\right| \leq \frac{1}{3} (T_2 - T_1)
\]
for every $T_1 \geq 1$ and $T_2 \geq T_1+1$ provided that $h \leq h_0$. On the other hand, by  Proposition \ref{spadaro0} (ii) we have  for every $h \leq h_0$ and $l \in \N$ 
\[
I_{l,h}:=\fint_{l^2}^{(l+1)^2} \|H_{E^h_t}-\lambda^h_t\|_{L^2(\pa E^h_t)}^2 \ \d t \leq \frac{C}{l}.
\]
By Chebysev's inequality
\[
\left|\{ t \in (l^2,(l+1)^2): \|H_{E^h_t}-\lambda^h_t\|_{L^2(\pa E^h_t)}^2 \geq 3 I_{l,h} \}\right| \leq \frac{1}{3} ((l+1)^2-l^2).
\]
Therefore,  by choosing $T_1=l^2$ and $T_2=(l+1)^2$  we deduce that the set 
\[
\left\{ t \in (T_1,T_2): |E^h_t| = \omega_{n+1}, \|H_{E^h_t}-\lambda^h_t\|_{L^2(\pa E^h_t)}^2 < 3 I_{l,h} \right\}
\]
is non-empty. Thus, if $h \leq h_0$, then there is a sequence of times $(T^h_l)_l$, with $l^2 \leq T^h_l \leq (l+1)^2$, such that  the corresponding sets satisfy $|E^h_{T^h_l}|=\omega_{n+1}$
and 
\beq
\label{L2decay}
\|H_{E^h_{T^h_l}}-\lambda^h_{T^h_l}\|_{L^2(\pa E^h_{T^h_l})} \leq C l^{-\frac12}.
\eeq

By slight abuse of the notation we set   $E^h_l := E^h_{T^h_l}$ and $\lambda_{l,h} :=\lambda^h_{T^h_l} $ for $h \leq h_0$. 
Since the sets $E^h_l$ are $C^2$-regular and bounded, then thanks to $P(E_0) \leq C_0$, $|E^h_l| \geq 1/C_0$, \eqref{L2decay} and Theorem \ref{thm1}
we find $l_0 = l_0(C_0,n)$ such that for every $l \geq l_0$ we have $1/C \leq \lambda_{l,h} \leq C$,
\beq
\label{thm 1.1 step 1 1}
|P(E_l^h) -  N_l^h(n+1)\omega_{n+1}(r_l^h)^n| \leq C l^{-\frac{q}{2}} \qquad \text{and} \qquad \sup_{E_l^h \Delta F_l^h} d_{\pa F_l^h} \leq C  l^{-\frac{q}{2}},
\eeq
where $r^h_l = n / \lambda_{l,h}$ and $F_l^h$ is a union of $N_l^h$-many pairwise disjoint (open) balls of radius $r_{l,h}$. 
Since $1/C \leq \lambda_{l,h} \leq C$, then also $1/C \leq r_{l,h} \leq C$, which together with  the perimeter estimate 
 $P(E_l^h)\leq P(E_0) \leq C_0$ implies that there is $N_0=N_0(C_0,n) \in \N$ such that $N^h_l \leq N_0$. Further
the distance estimate  in \eqref{thm 1.1 step 1 1},  together with  $1/C \leq r_{l,h}\leq C$ and $N^h_l \leq N_0$, yields
\[
|E_l^h \Delta F_l^h| \leq C  l^{-\frac{q}{2}}.
\]

Since $|E_l^h|= \omega_{n+1}$, then the  estimate above implies 
$|(r_{l,h})^{n+1}N_l^h-1| \leq C l^{-\frac{q}{2}}$ and further
$
| (r_{l,h})^n(N^h_l)^{\frac{n}{n+1}}-1| \leq C l^{-\frac{q}{2}}.
$
This inequality, the perimeter estimate in  \eqref{thm 1.1 step 1 1} and $N^h_l \leq N_0$   imply
\beq
\label{P-N}
|P(E_l^h) - (n+1)\omega_{n+1} (N_l^h)^\frac1{n+1}| \leq C l^{-\frac{q}{2}}.
\eeq
Since by Proposition \eqref{spadaro0} (ii) $(P(E^h_l))_{l \geq l_0}$ is non-increasing, then \eqref{P-N} implies that there is a positive integer $l_1= l_1(C_0,n) \geq l_0$ for which
$(N_l^h)_{l\geq l_1}$ is non-increasing for all $h \leq h_0$. 

\medskip

\textbf{Step 2:} \quad For $l \geq l_1$ and $h \leq h_0$ the sets  $E^h_l$  are thus close to $N_l^h$ many balls. We claim that there are $N \in \N$ and $l_2  \geq l_1$ such that for every integer $L \geq l_2$ it holds
\beq
\label{N sama}
N_l^h = N \qquad \text{for all } \, l_2 \leq l \leq L 
\eeq
provided that $h $ is small enough. 

Recall that by the assumption  $(E^h_t)_{t\geq 0}$ converges  to the flat flow $(E_t)_{t\geq 0}$
in $L^1$-sense compactly uniformly in time.  By using a standard diagonal argument  and possibly passing to a subsequence we 
find  a sequence of positive integers $(N_l)_{l\geq l_1}$, with $N_l \leq N_0$, such that 
$N^{h}_l \rightarrow N_l$ for every $l\geq l_1$. Since $(N^{h}_l)_{l \geq l_1}$ is non-increasing, then $(N_l)_{l\geq l_1}$ is non-increasing too and hence there are $N, l_2 \in \mathbb N$, $l_2 \geq l_1$, such that
$N_l = N$ for every $l \geq l_2$. Hence, we have \eqref{N sama} by the convergence of $N^{h}_l$ to $N_l$. 

We obtain from   \eqref{P-N} and \eqref{N sama} that 
\beq
\label{peri sama}
|P(E_l^h) - (n+1)\omega_{n+1} (N)^\frac1{n+1}| \leq C l^{-\frac{q}{2}}
\eeq
for $l_2 \leq l \leq L $, provided that  $h$ is small enough. Therefore, it follows from  Proposition \ref{spadaro0} (ii) that
\[
\int_{T^{h}_l + h}^{ T^{h}_L} \|H_{E^{h}_t} - \lambda^{h}_t\|_{L^2(\pa E^{h}_t)}^2 \ \d t \leq C  l^{-\frac{q}{2}}. 
\]
Since $h \leq 1$, and $L>1$ was arbitrary, the above yields  
\beq
\label{limsup}
 \sup_{T \geq (l+2)^2} \left[\limsup_{h \to 0} \int_{(l+2)^2}^T \|H_{E^{h}_t} - \lambda^{h}_t\|_{L^2(\pa E^{h}_t)}^2 \ \d t \right] \leq C  l^{-\frac{q}{2}}
\eeq
for every $l\geq l_2$.

\medskip

\textbf{Step 3:} \quad Let us fix small $\delta$, which choice will be clear later.  Then it follows from \eqref{limsup}, \eqref{peri sama} and the fact $t \mapsto P(E^h_t)$ is non-increasing in $\Sigma_h$ that there is $T_\delta$ such that for every $T\geq T_\delta + 1$ there is $h_{\delta, T}$ such that
\beq
\label{step 3 1} 
 \int_{T_\delta}^T \|H_{E^{h}_t} - \lambda^{h}_t\|_{L^2(\pa E^{h}_t)}^2 \ \d t  \leq \delta
\eeq
for all $h \leq h_{\delta, T}$ and
\beq \label{step 3 2}
|P(E_t^{h}) - (n+1)\omega_{n+1} N^\frac1{n+1}| \leq \delta 
\eeq
for all $t \in (T_\delta, T) \setminus \Sigma_h$.   On the other hand, by Proposition \ref{lambda bound} and by decreasing $ h_{\delta, T}$ if necessary  we deduce that
\beq \label{step 3 3}
| \Sigma_h \cap (T_\delta, T)| \leq \delta \qquad \text{for all } \, h \leq h_{\delta,T}.
\eeq

Let  $\eps>0$ and let us fix $t \geq T_\delta+1$. (The time $T_\delta+1$ will be  $T_\eps$  in  the claim.)   We claim  that, when $\delta$ is chosen small enough,   it holds
\beq 
\label{thm1.1 step 3} 
\sup_{E^{h}_t \Delta F_t^h} d_{\pa F_t^h} \leq \eps,
\eeq
for $h \leq h_{\delta, T}$, where $F_t^h$ is a union of $N$-many pairwise disjoint (open) balls of radius $r=N^{- \frac{1}{n+1}}$ with volume $\omega_{n+1}$.  

Fix $T \geq t+1$. Then it follows from \eqref{step 3 1}  that 
\[
 \int_{t-\delta^{1/4}}^t \|H_{E^{h}_\tau} - \lambda^{h}_\tau\|_{L^2(\pa E^{h}_\tau)}^2 \ \d \tau  \leq \delta 
\] 
and from \eqref{step 3 2} and \eqref{step 3 3} that 
\[
|P(E_\tau^{h}) - (n+1)\omega_{n+1} N^\frac1{n+1}| \leq \delta  \qquad \text{for all } \, \tau \in (t-\delta^{1/4}, t) \setminus \Sigma_h
\]
and $|\Sigma_h \cap (t-\delta^{1/4}, t)| \leq \delta $. Using these estimates  we deduce that there is  $t_0  \in (t-\delta^{1/4}, t)$ such that $|E_{t_0}^h| = \omega_{n+1}$, 
\beq 
\label{step 3 4} 
\big{|}P(E_{t_0}^{h}) - (n+1)\omega_{n+1} N^\frac1{n+1}\big{|}  \leq \delta
\eeq
and 
\[
\|H_{E^{h}_{t_0}} - \lambda^{h}_{t_0}\|_{L^2(\pa E^{h}_{t_0})} \leq \delta^{1/4}.
\]
Theorem \ref{thm1} implies that
\[
 \sup_{E^{h}_{t_0} \Delta F_{t_0}^h}   d_{\pa F_{t_0}^h} \leq C \delta^{q/4},
\]
for all $h \leq h_{\delta, T}$, where $F_{t_0}^h$  is a union of $N_{t_0, h}$-many pairwise disjoint (open) balls of radius $r_{t_0, h}$  with volume $\omega_{n+1}$ and 
\[
\big{|}P(E_{t_0}^{h}) - N_{t_0, h}(n+1)\omega_{n+1} r_{t_0, h}^n \big{|}  \leq C \delta^{q/4}.
\] 
ince $1/C \leq r_{t_0, h} \leq C$, then, as in Step 1, we deduce from the previous two estimates above that $|E^{h}_{t_0} \Delta F_{t_0}^h| \leq C \delta^{q/4}$.
Then by \eqref{step 3 4} and $|F_{t_0}^h| = \omega_{n+1}$  we further conclude that $N_{t_0, h} = N$, i.e.,  $F_{t_0}^h$  is a union of $N$-many pairwise disjoint (open) balls with volume $\omega_{n+1}$ and radius $r = N^{- \frac{1}{n+1}}$. 

By Lemma \ref{comparison} it holds 
\[
 \sup_{E^{h}_{\tau} \Delta F_{t_0}^h}   d_{\pa F_{t_0}^h} \leq C \delta^{\frac{q}{36}} \qquad \text{for all } t_0 < \tau < t_0 +  \delta^{\frac{q}{8}} 
\]
and $h \leq h_{\delta, T}$.  In particular, since $ \delta^{\frac{q}{8}} >  \delta^{\frac{1}{4}}$ the above inequality holds for $t$. This proves \eqref{thm1.1 step 3} by choosing $F_{t}^h = F_{t_0}^h$ and $\delta $  small enough. The claim follows by letting $h \to 0$. Note that by Proposition \ref{spadaro0} (iii) there is $R>0 $ such that $F_{t}^h \subset B_R$ for all $h \leq h_{\delta, T}$. Therefore, by passing to another subsequence if necessary, we have that  $F_{t}^h  \to F_{t} $, where $F_t$ is  a union of $N$-many pairwise disjoint (open) balls with volume $\omega_{n+1}$ and by  \eqref{thm1.1 step 3} it holds
\[
\sup_{E_t \Delta F_t} d_{\pa F_t} \leq \eps.
\]

\end{proof}


\section*{Acknowledgments}
The research was supported by the Academy of Finland grant 314227.


\end{document}